\magnification=\magstep1


\def\item{\vskip1.3pt\hang\textindent}


\tolerance=300
\pretolerance=200
\hfuzz=1pt
\vfuzz=1pt

\hoffset 0cm            
\hsize=5.8 true in
\vsize=9.5 true in

\def\rightheadline{\hfil\smc\lastname\hfil\tenbf\folio}
\def\leftheadline{\tenbf\folio\hfil\smc\lastname\hfil}
\headline={\ifodd\pageno\rightheadline\else\leftheadline\fi}
\newdimen\dimenone
\def\checkleftspace#1#2#3#4#5{
 \dimenone=\pagetotal
 \advance\dimenone by -\pageshrink   
 \ifdim\dimenone>\pagegoal          
   \else\dimenone=\pagetotal
        \advance\dimenone by \pagestretch
        \ifdim\dimenone<\pagegoal
          \dimenone=\pagetotal
          \advance\dimenone by#1         
          \setbox0=\vbox{#2\parskip=0pt                
                       \hyphenpenalty=10000
                       \rightskip=0pt plus 5em
                       \noindent#3 \vskip#4}    
        \advance\dimenone by\ht0
        \advance\dimenone by 3\baselineskip
        \ifdim\dimenone>\pagegoal\vfill\eject\fi
          \else\eject\fi\fi}

\parindent=35pt
\mathsurround=1pt
\parskip=1pt plus .25pt minus .25pt
\normallineskiplimit=.99pt

\mathchardef\emptyset="001F 

\def\Int{\mathop{\rm int}\nolimits}
%



\def\1{{\bf1}}\def\0{{\bf0}}

\def\({\bigl(}  \def\){\bigr)}
\def\<{\mathopen{\langle}}\def\>{\mathclose{\rangle}}

\def\Z{{\mathchoice{{\hbox{$\rm Z\hskip 0.26em\llap{\rm Z}$}}}%
{{\hbox{$\rm Z\hskip 0.26em\llap{\rm Z}$}}}%
{{\hbox{$\scriptstyle\rm Z\hskip 0.31em\llap{$\scriptstyle\rm Z$}$}}}{{%
\hbox{$\scriptscriptstyle\rm Z$\hskip0.18em\llap{$\scriptscriptstyle\rm Z$}}}}}}

\def\F{{\mathchoice{\hbox{$\rm I\hskip-0.14em F$}}%
{\hbox{$\rm I\hskip-0.14em F$}}%
{\hbox{$\scriptstyle\rm I\hskip-0.14em F$}}%
{\hbox{$\scriptscriptstyle\rm I\hskip-0.10em F$}}}}

\def\R{{\mathchoice{\hbox{$\rm I\hskip-0.14em R$}}%
{\hbox{$\rm I\hskip-0.14em R$}}%
{\hbox{$\scriptstyle\rm I\hskip-0.14em R$}}%
{\hbox{$\scriptscriptstyle\rm I\hskip-0.10em R$}}}}

\def\K{{\mathchoice{\hbox{$\rm I\hskip-0.15em K$}}%
{\hbox{$\rm I\hskip-0.15em K$}}%
{\hbox{$\scriptstyle\rm I\hskip-0.15em K$}}%
{\hbox{$\scriptscriptstyle\rm I\hskip-0.11em K$}}}}

\def\qed{\hfill {\hbox{[\hskip-0.05em ]}}}

\def\.{{\cdot}}
\def\|{\Vert}
\def\ssk{\smallskip}
\def\msk{\medskip}
\def\bsk{\bigskip}
\def\giantskip{\vskip2\bigskipamount}

\def\giantbreak{\par \ifdim\lastskip<2\bigskipamount \removelastskip
         \penalty-400 \giantskip\fi}

\def\nin{\noindent}
\def\cen{\centerline}
\def\pagebreak{\vskip 0pt plus 0.0001fil\break}
\def\linebreak{\break}

\def\epsilon{\varepsilon}

\font\ninerm=cmr9
\font\eightrm=cmr8
\font\sixrm=cmr6

\font\eightbf=cmbx8
\font\sixbf=cmbx6

\font\eighti=cmmi8
\font\sixi=cmmi6
\font\ninesy=cmsy9
\font\eightsy=cmsy8
\font\sixsy=cmsy6

\font\eightit=cmti8


\font\eightsl=cmsl8

\font\eighttt=cmtt8
\font\bfone=cmbx10 scaled\magstep1 
\font\smc=cmcsc10

\font\small=cmcsc8

\def\no #1. {\bigbreak\vskip-\parskip\noindent\bf #1. \quad\rm}

\def\Proposition #1. {\checkleftspace{0pt}{\bf}{Theorem}{0pt}{}
\bigbreak\vskip-\parskip\noindent{\bf Proposition #1.}
\quad\it}

\def\Theorem #1. {\checkleftspace{0pt}{\bf}{Theorem}{0pt}{}
\bigbreak\vskip-\parskip\noindent{\bf  Theorem #1.}
\quad\it}
\def\Corollary #1. {\checkleftspace{0pt}{\bf}{Theorem}{0pt}{}
\bigbreak\vskip-\parskip\nin{\bf Corollary #1.}
\quad\it}
\def\Lemma #1. {\checkleftspace{0pt}{\bf}{Theorem}{0pt}{}
\bigbreak\vskip-\parskip\noindent{\bf  Lemma #1.}\quad\it}

\def\Definition #1. {\checkleftspace{0pt}{\bf}{Theorem}{0pt}{}
\rm\bigbreak\vskip-\parskip\noindent{\bf Definition #1.}
\quad}

\def\Remark #1. {\checkleftspace{0pt}{\bf}{Theorem}{0pt}{}
\rm\bigbreak\vskip-\parskip\noindent{\bf Remark #1.}\quad}

\def\Exercise #1. {\checkleftspace{0pt}{\bf}{Theorem}{0pt}{}
\rm\bigbreak\vskip-\parskip\noindent{\bf Exercise #1.}
\quad}

\def\Example #1. {\checkleftspace{0pt}{\bf}{Theorem}{0pt}{}
\rm\bigbreak\vskip-\parskip\noindent{\bf Example #1.}\quad}
\def\Examples #1. {\checkleftspace{0pt}{\bf}{Theorem}{0pt}
\rm\bigbreak\vskip-\parskip\noindent{\bf Examples #1.}\quad}

\newcount\problemnumb \problemnumb=0
\def\Problem{\global\advance\problemnumb by 1\bigbreak\vskip-\parskip\noindent
{\bf Problem \the\problemnumb.}\quad\rm }

\def\Proof#1.{\rm\par\ifdim\lastskip<\bigskipamount\removelastskip\fi\smallskip
            \noindent {\bf Proof.}\quad}

\nopagenumbers

\def\author{}
\def\lastname{}
\def\thanks#1{\footnote*{\eightrm#1}}
\def\title{}

\def\lastname{}
\def\h{{\textstyle{1\over2}}}

\def\he{{1\over2}}

\def\n{{\cal N}}
\def\ep{\epsilon}

\def\text{\textstyle}
\def\disp{\displaystyle}
\def\d{{\,\rm d}}

\def\and{{\rm and }}
\def\und{{\rm und }}

\def\n{\cen{{\it W.G. Nowak}}}

\expandafter\edef\csname amssym.def\endcsname{%
       \catcode`\noexpand\@=\the\catcode`\@\space}
\catcode`\@=11
\def\undefine#1{\let#1\undefined}
\def\newsymbol#1#2#3#4#5{\let\next@\relax
 \ifnum#2=\@ne\let\next@\msafam@\else
 \ifnum#2=\tw@\let\next@\msbfam@\fi\fi
 \mathchardef#1="#3\next@#4#5}
\def\mathhexbox@#1#2#3{\relax
 \ifmmode\mathpalette{}{\m@th\mathchar"#1#2#3}%
 \else\leavevmode\hbox{$\m@th\mathchar"#1#2#3$}\fi}
\def\hexnumber@#1{\ifcase#1 0\or 1\or 2\or 3\or 4\or 5\or 6\or 7\or 8\or
 9\or A\or B\or C\or D\or E\or F\fi}

\font\tenmsb=msbm10
\font\sevenmsb=msbm7
\font\fivemsb=msbm5
\newfam\msbfam
\textfont\msbfam=\tenmsb
\scriptfont\msbfam=\sevenmsb
\scriptscriptfont\msbfam=\fivemsb
\edef\msbfam@{\hexnumber@\msbfam}
\def\Bbb#1{{\fam\msbfam\relax#1}}

\newsymbol\Bbbk 207C
\def\widehat#1{\setbox\z@\hbox{$\m@th#1$}%
 \ifdim\wd\z@>\tw@ em\mathaccent"0\msbfam@5B{#1}%
 \else\mathaccent"0362{#1}\fi}
\def\widetilde#1{\setbox\z@\hbox{$\m@th#1$}%
 \ifdim\wd\z@>\tw@ em\mathaccent"0\msbfam@5D{#1}%
 \else\mathaccent"0365{#1}\fi}
\font\teneufm=eufm10
\font\seveneufm=eufm7
\font\fiveeufm=eufm5
\newfam\eufmfam
\textfont\eufmfam=\teneufm
\scriptfont\eufmfam=\seveneufm
\scriptscriptfont\eufmfam=\fiveeufm

\catcode`@=11 

\expandafter\edef\csname amssym.def\endcsname{%
       \catcode`\noexpand\@=\the\catcode`\@\space}
\font\eightmsb=msbm8
\font\sixmsb=msbm6
\font\fivemsb=msbm5
\font\eighteufm=eufm8
\font\sixeufm=eufm6
\font\fiveeufm=eufm5
\newskip\ttglue
\def\eightpoint{\def\rm{\fam0\eightrm}%
  \textfont0=\eightrm \scriptfont0=\sixrm \scriptscriptfont0=\fiverm
  \textfont1=\eighti \scriptfont1=\sixi \scriptscriptfont1=\fivei
  \textfont2=\eightsy \scriptfont2=\sixsy \scriptscriptfont2=\fivesy
  \textfont3=\tenex \scriptfont3=\tenex \scriptscriptfont3=\tenex
\textfont\eufmfam=\eighteufm
\scriptfont\eufmfam=\sixeufm
\scriptscriptfont\eufmfam=\fiveeufm
\textfont\msbfam=\eightmsb
\scriptfont\msbfam=\sixmsb
\scriptscriptfont\msbfam=\fivemsb
  \def\it{\fam\itfam\eightit}%
  \textfont\itfam=\eightit
  \def\sl{\fam\slfam\eightsl}%
  \textfont\slfam=\eightsl
  \def\bf{\fam\bffam\eightbf}%
  \textfont\bffam=\eightbf \scriptfont\bffam=\sixbf
   \scriptscriptfont\bffam=\fivebf
  \def\tt{\fam\ttfam\eighttt}%
  \textfont\ttfam=\eighttt
  \tt \ttglue=.5em plus.25em minus.15em
  \normalbaselineskip=9pt
  \def\MF{{\manual opqr}\-{\manual stuq}}%
  \let\big=\eightbig
  \setbox\strutbox=\hbox{\vrule height7pt depth2pt width\z@}%
  \normalbaselines\rm}
\def\eightbig#1{{\hbox{$\textfont0=\ninerm\textfont2=\ninesy
  \left#1\vbox to6.5pt{}\right.\n@space$}}}


\csname amssym.def\endcsname


\def\la{\lambda}
\def\al{\alpha}
\def\be{\beta}

\def\({\left(}
\def\){\right)}

\def\eq{\eqalign}
\def\f{{1\over 2\pi i}}

\def\abs#1{\left| #1 \right|}

\def\norm#1{\left\Vert #1 \right\Vert}

\def\klein{\eightpoint \def\smc{\small} \baselineskip=9pt}

\def\fn#1#2{{\parindent=0.7true cm
\footnote{$^{(#1)}$}{{\klein  #2}}}}

\font\boldmas=msbm10                  
\def\Bbb#1{\hbox{\boldmas #1}}        
\def\Z{{\Bbb Z}}                        

\def\R{{\Bbb R}}
\def\F{{\Bbb F}}


\font\eightrm=cmr8
\long\def\fussnote#1#2{{\baselineskip=9pt
\setbox\strutbox=\hbox{\vrule height 7pt depth 2pt width 0pt}%
\eightrm
\footnote{#1}{#2}}}
\font\boldmasi=msbm10 scaled 700      
\def\Bbbi#1{\hbox{\boldmasi #1}}      
\font\boldmas=msbm10                  
\def\Bbb#1{\hbox{\boldmas #1}}        
\def\Zi{{\Bbbi Z}}                      
\def\Pi{{\Bbbi P}}                      
\def\Ri{{\Bbbi R}}



\def\dint #1 {
\quad  \setbox0=\hbox{$\disp\int\!\!\!\int$}
  \setbox1=\hbox{$\!\!\!_{#1}$}
  \vtop{\hsize=\wd1\centerline{\copy0}\copy1} \quad}

\def\drint #1 {
\qquad  \setbox0=\hbox{$\disp\int\!\!\!\int\!\!\!\int$}
  \setbox1=\hbox{$\!\!\!_{#1}$}
  \vtop{\hsize=\wd1\centerline{\copy0}\copy1}\qquad}

\def\frac#1#2{{#1\over #2}}

\def\date{\the\day.~\the\month.~\the\year}

\def\klein{\eightpoint \def\smc{\small} }

\def\frac#1#2{{#1\over#2}}
\def\Int{\int\limits}

\def\vol{{\rm vol}}


\hsize=16true cm     \vsize=24true cm

\parindent=0cm

\def\K{{\cal K}}
\def\B{{\cal B}}
\def\LP{11}  \def\AZ{12}
\def\b#1{{\bf #1}}
\def\n#1#2{\norm{#2}_{#1}}
\def\v{{4\pi\over3}}
\def\E{{\cal E}}
\def\M{{\cal M}}
\def\fou#1{\widehat{#1}}
\def\eb{^{(2)}} 

\vbox{\vskip1true cm}

\cen{\bfone Effektive Absch{\"a}tzungen f{\"u}r den Gitterrest}
\msk \cen{\bfone gewisser ebener und dreidimensionaler
Bereiche}\bsk

\cen{\bf Ekkehard Kr{\"a}tzel und Werner Georg Nowak (Wien)}
\bsk\bsk

\vbox{\vskip 1true cm}

\footnote{}{\klein{\it Mathematics Subject Classification }
(2000): 11P21, 11K38, 52C07.\par }

{\klein{\bf Abstract. \ Effective estimates for the lattice point
discrepancy of certain planar and three-dimensional domains. }
This paper provides estimates, with explicit constants, 
for the lattice point discrepancy of $\b o$-symmetric ellipse discs and
ellipsoids in $\Ri^3$, as well as of three-dimensional convex
bodies which are invariant under rotations around one coordinate
axis and have a smooth boundary of finite nonzero Gaussian
curvature throughout. }

\vbox{\vskip 1.2true cm}

{\bf 1.~Einleitung. } Das zentrale Problem der Theorie der
Gitterpunkte in gro{\ss}en Bereichen (wie sie z.B.~vom
erstgenannten Verfasser in der Monographie [\LP] dargestellt
wurde) besteht bekanntlich darin, f{\"u}r einen gegebenen
K{\"o}rper $\K$ des $\R^s$, \ $s\ge2$, und einen gro{\ss}en
reellen Parameter $t$ die Zahl $A_\K(t)$ der ganzzahligen Punkte
im linear vergr{\"o}{\ss}erten Bereich $t\,\K$ asymptotisch
auszuwerten, insbesondere den Gitterrest
$$ P_\K(t) = A_\K(t) - \vol(\K)t^s = \#\(\Z^s\cap t\,\K\) - \vol(\K)t^s  $$
m{\"o}glichst pr{\"a}zise abzusch{\"a}tzen. \ssk

Unter der Voraussetzung, da{\ss} der Rand $\partial\K$ von $\K$
hinreichend glatt und {\"u}berall von beschr{\"a}nkter, nicht
verschwindender Gau{\ss}scher Kr{\"u}mmung ist, bewiesen bereits
J.G.~Van der Corput [16] \ $P_\K(t)\ll t^{2/3}$ und E.~Hlawka [9] \
$P_K(t)\ll t^{s(s-1)/(s+1)}$ f{\"u}r $s\ge3$. Die sch{\"a}rfsten
bekannten Schranken lauten
$$ P_\K(t) \ll \ \cases{t^{131/208}(\log t)^{18637/8320} & f{\"u}r $s=2$,\cr
t^{63/43+\ep} & f{\"u}r $s=3$,\cr t^{40/17+\ep} & f{\"u}r $s=4$,
\cr t^{s-2+(s+4)/(s^2+s+2)+\ep} & f{\"u}r $s\ge5$.\cr}
\eqno(1.1)$$ Sie stammen von M.~Huxley [10] ($s=2$)
bzw.~W.~M{\"u}ller [14] ($s\ge3$). Im Spezialfall der Kugel ist
f{\"u}r $s\ge3$ Genaueres  bekannt, n\"amlich
$$ P(t) \ll\ \cases{t^{21/16+\ep} & f\"ur $s=3$,\cr t^2(\log t)^{2/3} &
f\"ur $s=4$,\cr t^{s-2} & f\"ur $s\ge5$.\cr} \eqno(1.2) $$ Diese
Absch\"atzungen stammen von D.R.~Heath-Brown [8] (als Verbesserung
fr\"uherer Ergebnisse von F.~Chamizo und H.~Iwaniec [4] sowie
von I.M.~Vinogradov [17]), bzw.~von A.~Walfisz [18]. F\"ur ein
allgemeines $\b o$-symmetrisches Ellipsoid $\E$ bewiesen Bentkus
und G\"otze [1]
$$ P_\E(t) \ll t^{s-2} \qquad \hbox{ f\"ur } s\ge9\,.
\eqno(1.3)$$ Dies wurde von F.~G\"otze [6] auf $s\ge5$ erweitert. 
F\"ur einen dreidimensionalen Rotationsk\"orper
(bez\"uglich einer Koordinatenachse) mit hinreichend glattem Rand
von {\"u}berall beschr{\"a}nkter, nicht verschwindender
Gau{\ss}scher Kr{\"u}mmung zeigte F.~Chamizo [3], sch\"arfer als
(1.1),
$$ P_\K(t) \ll t^{11/8}\,. \eqno(1.4)  $$ 

\bsk\msk

{\bf 2.~Ziele der vorliegenden Arbeit. } Die Verfasser wurden
anl{\"a}{\ss}lich von Vortr\"agen mehrfach gefragt, was (bei den
auch hier berichteten Resultaten) {\"u}ber die in den
$\ll$-Symbolen enthaltenen Konstanten ausgesagt werden k{\"o}nne.
Diese Problematik wurde von den meisten Experten der Theorie bis
dato str{\"a}flich vernachl{\"a}ssigt.\fn{1}{Als Ausnahme seien die bereits 
zitierten Arbeiten von Bentkus und G\"otze [1], [6] genannt. 
Die $\ll$-Konstante wird darin bis auf einen nur von $s$ 
abh\"angenden Faktor mittels der minimalen und maximalen 
Hauptkr\"ummungsradien des Ellipsoids explizit angegeben.} 
Aus dem Blickwinkel des
allgemein interessierten Mathematikers sowie des Numerikers
verdient Sie jedoch durchaus Beachtung. \ssk Dabei liegt es
allerdings nahe, sich mit effektiven Versionen der klassischen
Schranken von Van der Corput bzw.~Hlawka zu begn{\"u}gen. In
dieser Richtung wurde von E.~Kr{\"a}tzel [12], Satz 5.12, f{\"u}r den
Einheitskreis die Absch\"atzung 
$$ \abs{P(t)} \le 38\,t^{2/3}+704\,t^{1/2} + 11 \eqno(2.1) $$
gezeigt. Allgemeiner konnte er in [13] f{\"u}r einen ebenen Bereich
$\B$ mit glattem Rand von stetiger, beschr{\"a}nkter, nicht
verschwindender Kr{\"u}mmung
$$ \abs{P_\B(t)} \le 48\,\(r_{\max}\,t\)^{2/3} + \(703\,\sqrt{r_{\max}}
+{3\,r_{\max}\over5\,\sqrt{r_{\min}}}\)\,\sqrt{t} + 11 \eqno(2.2)
$$ beweisen, wobei $r_{\max}$, $r_{{\min}}$ die Extremwerte des
Kr{\"u}mmungsradius der Randkurve bezeichnen. \ssk

In dieser Arbeit soll die Thematik auf gewisse dreidimensionale
K\"orper erweitert werden. Wir ben\"utzen zun\"achst (Abschnitt 3)
die Poisson'sche Summenformel zusammen mit einer
Gl\"attungstechnik und Kenntnissen \"uber die auftretenden
Fourier-Transformierten, um f\"ur den Gitterrest der Kugel sowie
eines allgemeinen $\b o$-symmetrischen Ellipsoids im $\R^3$
Schranken mit effektiven Konstanten zu erhalten (S\"atze 1A, 1B).
Dieser Zugang liefert auch ein entsprechendes Resultat f\"ur die
ebene Ellipsenscheibe (Satz 2). Mit ganz anderer Methodik,
n\"amlich einer Anwendung des Van der Corput'schen Ideenkreises
mittels Exponentialsummen und Bruchteilsummen (unter Einbeziehung
der Kusmin-Landau'schen Absch\"atzung), wird sodann im
4.~Abschnitt der Gitterrest eines dreidimensionalen
Rotationsk\"orpers mit glattem Rand von beschr\"ankter, nicht
verschwindender Gau\ss scher Kr\"ummung behandelt (Satz 3).
\bsk\bsk

\parindent=0cm

\def\b#1{{\bf #1}}
\def\n#1#2{\norm{#2}_{#1}}
\def\v{{4\pi\over3}}
\def\E{{\cal E}}
\def\M{{\cal M}}
\def\fou#1{\widehat{#1}}
\def\eb{^{(2)}}

\vbox{{\bf 3.~Gitterpunkte in Kugeln, Ellipsoiden und
Ellipsenscheiben.} \bsk

{\bf Satz 1A.}\quad{\it Es sei
$$  A(t) = \#\{ \b m \in \Z^3:\ \n2{\b m} \le t\ \}\,,\quad P(t) =
A(t)- \v t^3\,,  $$ dann gilt f\"ur $t\ge10$ die Absch\"atzung
$$ \abs{P(t)} \le 14 t^{3/2} + 3.5\, t \log t + 12.9\, t +
0.53\sqrt{t} + 8.2 + {2\log t\over\sqrt{t}} + {7.3\over\sqrt{t}} +
{0.36\over t}+ {1\over t^{3/2}}\,. $$ }} \bsk

Es sei dem Leser \"uberlassen, die gegebene Schranke auf Kosten
der Pr\"azision beliebig zu vereinfachen.

Unser Beweis l\"a\ss t sich prinzipiell in gr\"o\ss erer
Allgemeinheit durchf\"uhren, soda\ss\ wir den Fall eines
beliebigen dreidimensionalen Ellipsoids (in Mittelpunktslage)
behandeln k\"onnen.

F\"ur $\b x \in \R^3$  sei $ Q = Q(\b x) = \b x A \,^t\b x $ eine
tern\"are positiv definite quadratische Form mit Determinante
$\det A = 1$, und $f=\sqrt{Q}$ die Distanzfunktion des Ellipsoids
$\E=\E_f:\ Q(\b x)\le 1$.

Wir setzen $b:=
{1\over\pi}\sqrt{{3\over8}+{\sqrt{15}\over8}}=0.295\dots$,
$c:=0.277$, und definieren f\"ur $t>0$
$$ \ep := {c\over\sqrt{t}}\,,\quad K := {b\over\ep}\,. \eqno(3.1) $$
Ferner sei $Q^{-1}(\b x)=\b x A^{-1} \,^t\b x $ die inverse Form
zu $Q$ und entsprechend $g=\sqrt{Q^{-1}}$ die {\it
St\"utzfunktion} von $\E$. Wir setzen
$$ G_0 := \min_{\b o\ne\b k\in\Zi^3} g(\b k)\,, \quad G_1 :=
\max_{\b w\in[-\he,\he]^3} g(\b w) = \max g(\pm\h,\pm\h,\pm\h)\,,
\eqno(3.2)$$ und schlie\ss lich, f\"ur reelles $u\ge0$,
$$ a_g(u) := \v(u+G_1)^3 -1\,. \eqno(3.3)  $$ \bsk

{\bf Satz 1B.}\quad{\it Mit den vereinbarten Voraussetzungen und
Bezeichnungen sei
$$ A_f(t) := \#\(t\E \cap \Z^3\)\,,\quad P_f(t):=A_f(t) - \v
t^3\,, $$ dann gilt f\"ur $t\ge G_0^2$
$$ \abs{P_f(t)} \le \v\((t+2\ep)^3-t^3\) 
+ \max_{\pm}\,\overline{P_f^*}(t\pm2\ep) $$ mit
$$  \overline{P_f^*}(R) := {R\over\pi}\,{a_g(K)\over K^2} +
{2R\over\pi}\Int_{G_0}^K {a_g(u)\over u^3} \d u + {1\over2\pi^3 R}
\Int_{G_0}^\infty {a_g(u)\over u^5} \d u + \sum_{j=6,8,10,12}
F_j\, j\Int_K^\infty{a_g(u)\over u^{j+1}} \d u\,,  $$ wobei
$R:=t\pm2\ep$ und
$$ F_6 := {9R\over16\pi^5 \ep^4}\,,\quad F_8 := {9(\ep^2+2R^2)\over 128\pi^7 \ep^6
R}\,,\quad F_{10} := {9(2\ep^2+R^2)\over1024\pi^9\ep^8 R}\,,\quad
F_{12} := {9\over8192\pi^{11}\ep^8 R }\,.$$ } \bsk

{\bf Bemerkung.}\quad Diese Absch\"atzung erscheint auf den ersten
Blick reichlich kompliziert und nicht-explizit. Tats\"achlich ist
es aber f\"ur jede konkrete Form $Q$ leicht, $G_0$ und $G_1$ zu
bestimmen. Dann h\"angt die gegebene Schranke nur mehr von $t$ ab,
die Integrale k\"onnen (z.B.~mit Unterst\"utzung von {\it Derive}
[15]) elementar ausgewertet werden, und man erh\"alt stets ein
\"ahnliches Ergebnis wie in Satz 1A, wobei der gr\"o\ss te Term
durchwegs $14 t^{3/2}$ ist unabh\"angig von $Q$. (Vgl.~Bemerkung 2
am Ende des Beweises.)\ssk Um Satz 1A aus Satz 1B herzuleiten,
braucht man letzteren nur auf den Spezialfall der Kugel
anzuwenden. Hier ist $G_0=1,\ G_1=\h\sqrt{3}$, und man erh\"alt
durch Auswerten der Integrale und Einsetzen aller Definitionen die
in Satz 1A behauptete Schranke, abgesehen von etwas besseren
numerischen Konstanten und einem zus\"atzlichen Summanden
$$ h(t) = {0.6302963\,t+1.9074383\sqrt{t}+0.23678722\over6.878339\, t^{3/2}-3.8106}\,.
$$ Dieser f\"allt monoton in $t$, und es ist $h(t)\le0.06$ f\"ur $t\ge10$.
Daraus folgt Satz 1A trivial. \bsk\msk

{\bf Beweis von Satz 1B.}\quad Wir verwenden fettgedruckte
Kleinbuchstaben f\"ur Elemente des $\R^3$, insbesondere seien $\b
m, \b k \in \Z^3 $. F\"ur beliebige Teilmengen $\M \subseteq\R^3$
bezeichne $\chi_\M$ die Indikatorfunktion. Weiters sei
$$  \delta_0(\b x) := \cases{{3\over4\pi}\,{1\over\ep^3} &wenn $f(\b x) \le\ep\,,$\cr 0 &
sonst,\cr}$$
$$ \delta(\b x) := (\delta_0*\delta_0)(\b x) = \Int_{\Ri^3}
\delta_0(\b y) \delta_0(\b x-\b y)\d \b y = \Int_{\Ri^3}
\delta_0(\b x +\b y) \delta_0(\b y)\d \b y\,.  $$ Unsere Strategie
ist, die Indikatorfunktion des Ellipsoids $t\E$ durch Faltung mit
$\delta$ zu einer stetigen Funktion zu gl\"atten, soda\ss\ darauf
die Poissonsche Formel angewendet werden kann und zu einer
konvergenten Reihe f\"uhrt. \msk

\vbox{Wir behaupten zun\"achst: \msk

\qquad(i) Aus $f(\b x)>2\ep$ folgt stets $\delta(\b x)=0$. \ssk

\qquad(ii) Es gilt \quad $\chi_{(t-2\ep)\E}\,*\delta \le
\chi_{t\E}\le \chi_{(t+2\ep)\E}\,*\delta $.} \msk

{\it Beweis.}\quad (i) Es sei $f(\b x)>2\ep$, $\delta_0(\b y) \ne
0$, dann folgt $f(\b y)\le \ep$ und weiter, nach der
Dreiecksungleichung f\"ur die konvexe Funktion $f$, $f(\b x +\b
y)\ge f(\b x)-f(\b y)>\ep$, daher $\delta_0(\b x +\b y) = 0$. Also
ist $\delta_0(\b x +\b y) \delta_0(\b y)=0$ identisch in $\b y$,
damit (i) gezeigt. \ssk (ii) 1.~Fall: Es sei $\chi_{t\E}(\b x)=0$,
also $f(\b x)>t$. Ist au\ss erdem $\delta(\b y)\ne0$, dann folgt
$f(\b y)\le2\ep$ nach (i), und damit $f(\b x+\b y)>t-2\ep$. Daher
ist $\chi_{(t-2\ep)\E}(\b x +\b y)\delta(\b y)=0$ identisch in $\b
y$, folglich $(\chi_{(t-2\ep)\E}\,*\delta)(\b x)=0$. Der rechte
Teil von (ii) ist trivial. \ssk 2.~Fall: Es sei $\chi_{t\E}(\b
x)=1$, also $f(\b x)\le t$. Ist au\ss erdem $\delta(\b y)\ne0$,
dann folgt wieder $f(\b y)\le2\ep$, und weiter $f(\b x+\b y)\le
t+2\ep$. Daher ist $\chi_{\Ri^3\setminus(t+2\ep)\E}(\b x +\b
y)\delta(\b y)=0$ identisch in $\b y$, folglich
$(\chi_{\Ri^3\setminus(t+2\ep)\E}\,*\delta)(\b x)=0$, also
$(\chi_{(t+2\ep)\E}\,*\delta)(\b x)=1$. Der linke Teil ist jetzt
trivial. \qed \msk

Summation \"uber alle $\b m \in \Z^3$ ergibt
$$ \sum_{\b m\in \Zi^3}(\chi_{(t-2\ep)\E}\,*\delta)(\b m) \le
A_f(t) \le \sum_{\b m\in \Zi^3} (\chi_{(t+2\ep)\E}\,*\delta)(\b
m)\,. $$ Mit der mehrdimensionalen Poissonschen Summenformel
(vgl.~Bochner [2]) erhalten wir daraus
$$ \sum_{\b k\in \Zi^3}\fou{\chi_{(t-2\ep)\E}}(\b k)\,
(\fou{\delta_0}(\b k))^2 \le
A_f(t) \le \sum_{\b k\in \Zi^3} \fou{\chi_{(t+2\ep)\E}}(\b k)
\,(\fou{\delta_0}(\b k))^2\,, \eqno(3.4) $$  wenn \ $\fou{\cdot}$ \ wie
\"ublich die Fouriertransformierte bedeutet. Um letztere
auszuwerten, bemerken wir zun\"achst, da\ss\ $\fou{\delta_0}
(\b o)=1$ gilt und $\fou{\chi_{(t\pm2\ep)\E}}(\b o)=\v(t\pm2\ep)^3$.
Setzen wir (mit $e(w):=e^{2\pi iw}$)
$$ I_f(\b z) := \Int_{f(\b u)\le1} e(\b z \cdot\b u) \d\b u\,,  $$
dann ist weiters, f\"ur $\b k\ne\b o$ und $R=t\pm2\ep$,
$$  \fou{\chi_{R \E}}(\b k)= R^3\,I_f(R\b k)\,,\quad
\fou{\delta_0}(\b k) = {3\over4\pi}\,I_f(\ep\b k)\,. \eqno(3.5)
$$ Daraus folgt
$$ \abs{P_f(t)} \le \v\((t+2\ep)^3-t^3\) 
+ \max_{\pm}\,\abs{P_f^* (t\pm2\ep)} \eqno(3.6)  $$ mit
$$  P_f^*(R) := R^3\,\sum_{\b o\ne\b k\in \Zi^3}
I_f(R\b k)\,\({3\over4\pi}\,I_f(\ep\b k)\)^2\,.  \eqno(3.7) $$ Es
bleibt $I_f(\b z)$ f\"ur $\b z\ne\b o$ abzusch\"atzen. Trivial
gilt
$$ \abs{I_f(\b z)} \le \v\,.  \eqno(3.8) $$
Sei andererseits $B$ eine (symmetrische) \ $(3\times3)$-Matrix mit
$B^2=A$, dann gilt $$f^2(\b u)=Q(\b u)=\b u\, A\,^t\b u = (\b u
B)\,^t(\b u B)\,.$$ Daher ergibt die lineare Transformation $\b
v=\b u B$, wenn wir $\b w := \b z B^{-1}$ setzen, wegen $\b
z\cdot\b u = \b z\, ^t(\b v B^{-1}) = \b w \,^t\b v=\b w\cdot\b v
$,
$$ I_f(\b z) = \Int_{\n2{\b v}\le1} e(\b w\cdot\b v)\d\b v\,. $$
Drehen wir noch das Koordinatensystem so, da\ss\ die neue
$v_1$-Achse den Vektor $\b w$ enth\"alt, dann folgt
$$ I_f(\b z) = \Int_{\n2{\b v}\le1} e(\n2{\b
w}\,v_1)\d(v_1,v_2,v_3) = \Int_{\n2{\b v}\le1} e(g(\b z
)\,v_1)\d(v_1,v_2,v_3)\,,$$ wegen $\n2{\b w}^2 = \b w\, ^t\b w =
\b z B^{-1}\,^t(\b z B^{-1})= \b z A^{-1}\,^t\b z=g^2(\b z)$. Wir
schlie\ss en weiter, da\ss  \def\g{g(\b z)}
$$ \abs{I_f(\b z)} = \pi\,\abs{\Int_{-1}^1 (1-v_1^2)\,e(g(\b z) v_1) \d
v_1} = \abs{{\sin(2\pi g(\b z))\over2\pi^2 g(\b z)^3}-{\cos(2\pi
g(\b z) )\over\pi g(\b z)^2}} \le $$
$$ \le {1\over\pi\g^2}\,\sqrt{1+{1\over4\pi^2\g^2}}\ \le\
{1\over\pi\g^2}+ {1\over8\pi^3\g^4}\,.  \eqno(3.9) $$ Verwenden
wir dies in (3.7), so folgt
$$  \eq{\abs{P_f^*(R)} \le R^3\,\sum_{\b o\ne\b k\in \Zi^3} &
\({1\over\pi R^2 g(\b k)^2}+{1\over8\pi^3 R^4 g(\b k)^4}\)\times
\cr & \times\(\min\(1,{3\over4\pi^2\ep^2 g(\b k)^2} +
{3\over32\pi^4\ep^4 g(\b k)^4}\)\)^2\,.\cr} \eqno(3.10)  $$ Nach
Definition von $K$ und $b$ (vgl.~(3.1)) ist das auftretende
Minimum gleich 1 genau f\"ur $g(\b k)\le K$. Wir teilen daher die
Summe wie folgt auf:
$$  \abs{P_f^*(R)} \le R^3\,\sum_{0<g(\b k)\le K} \dots \
+ \ R^3 \sum_{g(\b k)>K} \dots\ =: S_1 +S_2\,.  \eqno(3.11)$$
Setzen wir $A_g^-(u) := A_g(u)-1$, dann ergibt sich mit Verwendung
von Stieltjes-Integralen (wegen $t\ge G_0^2$ ist sicher $K>G_0$)
\def\au{A_g^-(u)}
$$ S_1 \le {R\over\pi} \Int_{G_0-}^K u^{-2}\d A_g^-(u) + {1\over8\pi^3 R}
\Int_{G_0-}^\infty u^{-4}\d\au = $$
$$ = {R\over\pi}\,{A_g^-(K)\over K^2} +
{2R\over\pi}\Int_{G_0}^K {\au\over u^3} \d u + {1\over2\pi^3 R}
\Int_{G_0}^\infty {\au\over u^5} \d u\,. \eqno(3.12)$$ Definieren
wir, f\"ur $j\ge4$,
$$ s_j(g,t) := \sum_{g(\b k)>K} g(\b k)^{-j}\,,  $$
dann folgt einerseits durch partielle Integration
$$ s_j(g,t) = \Int_{K+}^\infty u^{-j}\,\d\au
\le j \Int_K^\infty {\au\over u^{j+1}}\d u\,, \eqno(3.13)$$
andererseits, mit R\"uckblick auf (3.10), (3.11) und die
Definition der $F_j$ in Satz 1B,
$$  S_2 = \sum_{j=6,8,10,12} F_j\,s_j(g,t)\,,  \eqno(3.14)  $$
wieder vorzugsweise mit Unterst\"utzung von {\it Derive } [15]. Zum
Beweis von Satz 1B ist es nun nur noch erforderlich, die
Absch\"atzung
$$\au\le a_g(u)= \v(u+G_1)^3 -1 \eqno(3.15)$$
(vgl.~(3.3)) zu verifizieren. Dazu zeigen wir, da\ss
$$ \bigcup_{\b k\in\Zi^3,\ 0\le g(\b k)\le u} (\b k + [-\h,\h]^3)
\subseteq (u+G_1)\E_g \eqno(3.16)  $$ gilt. Ist n\"amlich $\b
x=\b k+\b w$ mit $g(\b k)\le u$, $\n\infty{\b w}\le\h$, dann folgt
wegen (3.2) $g(\b w)\le G_1$ und wegen der Konvexit\"at $g(\b k+\b
w)\le u+G_1$, also $\b x\in(u+G_1)\E_g$. Betrachtung der Volumina
in (3.16) zeigt sofort (3.15). Zusammen mit (3.6) und (3.11)
bis (3.14) ergibt sich unmittelbar Satz 1B. \msk {\bf
Bemerkungen.}\quad1. Man k\"onnte die erzielte Absch\"atzung
weiter versch\"arfen, indem man anstelle der groben Schranke
(3.15) unseren Satz 1B iterativ in das Argument einsetzt. Auf
diese Weise w\"urde z.B.~der Term mit $t\log t$ in Satz 1A
verschwinden, allerdings bliebe der Hauptfehlerterm $14\, t^{3/2}$
unver\"andert. Wir begn\"ugen uns damit, f\"ur diesen gr\"o\ss ten
Term einen m\"oglichst kleinen Koeffizienten erhalten zu haben,
dazu einige weitere Summanden von kleinerer Ordnung, die nicht
wesentlich st\"oren.\ssk 2. Es ist instruktiv, sich im Detail
anzusehen, wie der erw\"ahnte Hauptfehler $14\, t^{3/2}$
unabh\"angig von der quadratischen Form $Q$ zustande kommt, und
warum die Wahl $c=0.277$ vor (3.1) ebenfalls f\"ur jedes Ellipsoid
in diesem Sinn optimal ist. Wir verwenden dazu die Notation (f\"ur
beliebige reelle Funktionen $H_1$ und $H_2>0$)
$$ H_1(t) \prec H_2(t) \quad :\iff\quad \lim_{t\to\infty}
{\max(H_1(t)-H_2(t),0)\over H_2(t)}= 0\,.  $$ So gilt nach
Konstruktion $R^{\pm1} \prec t^{\pm1}$, $a_g(u)\prec \v u^3$,
damit folgt aus der Darstellung in Satz 1B nach leichter Rechnung
$$  \overline{P_f^*}(t\pm2\ep) \prec
\(4 b +{3\over2\pi^4 b^3 }+{3\over10\pi^6 b^5} +{15\over896\pi^8
b^7}\){t^{3/2}\over c}\,. $$ Au\ss erdem gilt nat\"urlich $
\v\((t+2\ep)^3-t^3\)\prec 8\pi c t^{3/2} $, also
$$ P_f(t) \prec 8\pi c t^{3/2} + \(4 b +{3\over2\pi^4 b^3 } +
{3\over10\pi^6 b^5} + {15\over896\pi^8 b^7}\){t^{3/2}\over c}\,.
\eqno(3.17) $$ Ausbalanzieren bez\"uglich $c$ und Einsetzen des
Wertes von $b$ (vor (3.1)) f\"uhrt auf $c\approx0.277$, und die
rechte Seite von (3.17) wird $13.92\dots t^{3/2}$.\bsk

Die hier verwendete Methode kann auch verwendet werden, um die
eingangs zitierte Absch\"atzung (2.1) etwas zu verbessern und zu
verallgemeinern. Wir formulieren dieses Ergebnis als \msk

{\bf Satz 2.}\quad{\it Es sei $Q$ eine bin\"are positiv definite
quadratische Form mit Determinante 1, und wie vorher $f=\sqrt{Q}$
und $g=\sqrt{Q^{-1}}$ Distanzfunktion bzw.~St\"utzfunktion der
Ellipse $\E\eb:\ Q\le1$. Weiters sei analog fr\"uher
$$ A_f\eb(t) := \#\(t\E\eb \cap \Z^2\)\,,\quad P_f\eb(t):=A_f\eb(t) -
\pi t^2 $$  und  $$  G_0 := \min_{\b o\ne\b k\in\Zi^2} g(\b k)\,,
\quad G_1 := \max_{\b w\in[-\he,\he]^2} g(\b w) \,, \quad
a_g\eb(u) := \pi(u+G_1)^2-1\,. $$ Dann gilt f\"ur $t\ge 28\,
G_0^3$
$$ \eq{\abs{P_f\eb(t)} \le &\ \pi\((t+\ep)^2-t^2\) +
{(t+\ep)^{1/2}\over3K^{3/2}}\, a_g\eb(K) +\cr + &
{(t+\ep)^{1/2}\over2} \Int_{G_0}^K u^{-5/2} a_g\eb(u) \d u +
{(t+\ep)^{1/2}\over3\pi\ep^{3/2}} \Int_K^\infty u^{-4} a_g\eb(u)
\d u\,,\cr} \eqno(*)$$ wobei jetzt
$$ \ep:= {c\over t^{1/3}}\,,\quad c:={1\over6}\(2823576\over\pi^2\)^{1/9}
= 0.67317\dots\,, \quad  K = (3\pi)^{-2/3}\ep^{-1} $$ sei.
Insbesondere folgt f\"ur den Fall des Kreises, f\"ur $t\ge28$,
$$\abs{P\eb(t)}\le 8.46\, t^{2/3} + 1.5\, t^{1/2} + 1 + {3\over
t^{1/3}}+{2.85\over t^{2/3}} + {0.51\over t^{5/6}} + {0.34\over
t^{4/3}}\,.
$$ }\bsk
{\bf Beweisskizze zu Satz 2.} \quad Wir verwenden weitgehend das
fr\"uhere Argument und notieren nur die notwendigen
Modifikationen. Zun\"achst k\"onnen wir direkt f\"ur $\b x\in\R^2$
$$  \delta(\b x) := \cases{{1\over\pi}\,{1\over\ep^2} &wenn
$f(\b x) \le\ep\,,$\cr 0 &sonst,\cr}$$ setzen, eine Faltung mit
sich selbst ist nicht mehr n\"otig. Anstelle von (3.4) erh\"alt
man somit
$$ \sum_{\b k\in \Zi^2}\fou{\chi_{(t-\ep)\E\eb}}(\b k)\,\fou{\delta}(\b k) \le
A_f\eb(t) \le \sum_{\b k\in \Zi^2} \fou{\chi_{(t+\ep)\E\eb}}(\b k)
\,\fou{\delta}(\b k)\,. \eqno(3.18) $$ Dabei ist analog zu
(3.5), f\"ur $\b o\ne\b k\in\Z^2$ und $R=t\pm\ep$,
$$  \fou{\chi_{R \E\eb}}(\b k)= R^2\,I_f(R\b k)\,,\quad
\fou{\delta}(\b k) = {1\over\pi}\,I_f(\ep\b k)\,. \eqno(3.19)$$
Allerdings f\"uhrt die Auswertung des Integrals jetzt auf
$$ \abs{I_f(\b z)} = 2 \abs{\Int_{-1}^1 \sqrt{1-v^2}\,e(g(\b z)v)\d v}
= \abs{J_1(2\pi g(\b z))\over g(\b z)}\,, $$ wobei $J_1$ die
\"ubliche Besselfunktion bezeichnet. Nun ist\fn{2}{Bekanntlich
gilt $J_1(x)=(2/\pi)^{1/2}\,x^{-1/2}\cos(x-3\pi/4)+O(x^{-3/2})$,
wobei das Restglied z.B.~nach Gradshteyn \& Ryzhik [7], F.~8.451,
effektiv abgesch\"atzt werden kann. Das Maximum wird tats\"achlich
im ersten relativen Extremum $x_{\rm max}=2.16587\dots$
angenommen.}
$$ \max_{x>0} \abs{J_1(x)\sqrt{x}} = 0.82503\dots\,,  \eqno(3.20)$$
woraus wir
$$ \abs{I_f(\b z)} \le 0.83\,(2\pi)^{-1/2} g(\b z)^{-3/2}<{\text{1\over3}}
g(\b z)^{-3/2} $$ folgern. Mit (3.18), (3.19) ergibt sich daher
$$ \abs{P_f\eb(t)} \le \pi\((t+\ep)^2-t^2\) +
{(t+\ep)^{1/2}\over3\pi}\sum_{\b o\ne\b k\in\Zi^2 } g(\b
k)^{-3/2}\,\min\(\pi, {\text{1\over3}}\ep^{-3/2}g(\b
k)^{-3/2}\)\,. $$ Das Minimum in dieser Summe f\"uhrt auf deren
Aufteilung, je nachdem ob $g(\b k)\le K$ oder $g(\b k)> K$ ist.
Partielle Integration und eine elementare Absch\"atzung analog zu
(3.16) vervollst\"andigen den Beweis von Satz 2. \par Im
Spezialfall des Kreises ist nat\"urlich $G_0=1$,
$G_1={\sqrt{2}\over2}$. Mit der Ungleichung
$(t+\ep)^{1/2}\le\sqrt{t}\(1+{\ep\over2t}\)$ erh\"alt man die
Aussage von Satz 2 auch f\"ur diesen Fall.\ssk Um wieder das
Entstehen des Hauptfehlerterms zu \"uberblicken und die Wahl von
$\ep$ einzusehen, bemerken wir, da\ss\ $a_g\eb(u)\prec \pi u^2$
gilt. Damit erh\"alt man aus $(*)$
$$ P_f\eb(t) \prec 2\pi\ep t + {7\over9}(3\pi)^{2/3}\sqrt{t\over\ep}\,.
$$ Ausbalanzieren ergibt $\ep$ wie in Satz 2 und damit f\"ur den Hauptfehler
$P_f\eb(t)\prec 8.46\,t^{2/3}$, wieder unabh\"angig von der Form
$Q$.

\bsk\bsk

\hsize=16true cm     \vsize=23.5true cm

\parindent=0cm

\def\K{{\cal K}} \def\c{{\cal C}} \def\sgn{{\rm sgn}\,}
\def\LP{LP}  \def\AZ{AZ}
\def\ut{\left(\frac{u}{t}\right)}  \def\sumi{\mathop{{\sum}'}}

\def\f{{\textstyle{\,\min_{[a,b]}|f''|}}}
\def\F{{\textstyle{\,\max_{[a,b]}|f''|}}}
\def\ff{\({\textstyle{\min_{[a,b]}|f''|}}\)}

\def\rmax{r_{\max}} \def\rmin{r_{\min}}
\def\rq{\,{\rmax^3\over\rmin^4}\,}
\def\rf#1#2{\rmax^{#1}\,F_2^{#2}}

\def\iz#1{\Int_{z_j}^{z_{j+1}} |f''(w)|^{#1} }
\def\sj{\sum_{j=0}^J}


{\bf 4.~Gitterpunkte in einem dreidimensionalen
Rotationsk{\"o}rper.} \msk Wir untersuchen nun die analoge
Problemstellung (Absch{\"a}tzung des Gitterrestes mit expliziten
Konstanten) f{\"u}r einen kompakten konvexen K{\"o}rper $\K$ im
$\R^3$, der rotationssymmetrisch bez{\"u}glich der $z$-Achse und
au\ss erdem spiegelsymmetrisch bez{\"u}glich der $(x,y)$-Ebene
sei. Sein Rand $\partial\K$ sei hinreichend glatt und {\"u}berall
von beschr{\"a}nkter, nicht-verschwindender Gau\ss scher
Kr{\"u}mmung. Dies pr{\"a}zisieren wir am bequemsten, in dem wir
die Schnittkurve $\c$ von $\partial\K$ mit der $(y,z)$-Ebene
betrachten: Diese sei z.B.~parametrisiert durch
$(y,z)=(\rho(\theta)\sin\theta, \rho(\theta)\cos\theta)$, wobei
$\rho(\cdot)$ eine gerade, positive, periodische Funktion mit
Periode $\pi$ sei, {\"u}berall viermal stetig differenzierbar, und
durchwegs\fn{3}{Dies entspricht dem Nichtverschwinden der
Kr{\"u}mmung.} $\rho\,\rho''- 2\rho'^2-\rho^2\ne0$. Dann ist
$$ r=r_{\theta}=\abs{\rho(\theta)\,\rho''(\theta)-
2\rho'^2(\theta)-\rho^2(\theta)}^{-1}\,
(\rho^2(\theta)+\rho'^2(\theta))^{3/2} $$ der entsprechende
Kr{\"u}mmungsradius. Wir setzen o.b.d.A.~$\rho(0)=1$ voraus und
weiters\fn{4}{Diese zus{\"a}tzliche Annahme ist bequem und
vereinfacht die sp{\"a}teren ohnehin recht volumin{\"o}sen
Formeln. Unser Argument gilt aber unver{\"a}ndert f{\"u}r
$\rho({\pi\over2})>1$, wobei allerdings dann das Endergebnis auch
von $\rho({\pi\over2})$ abh{\"a}ngt.} $\rho({\pi\over2})\le1$. Die
Kurve $\c$ kann andererseits in den vier Quadranten durch
nicht-negative Funktionen\fn{5}{Es besteht wohl keine
Verwechslungsgefahr mit den fr{\"u}her so bezeichneten
Distanzfunktionen.} $f, g$ in der Form $y=\pm f(z)\iff z=\pm g(y)$
dargestellt werden. Taylorentwicklung um $y=0$ ergibt $z=g(y)=1+\h
r_0^{-1}\,y^2+O(\abs{y}^3)$ und damit, f{\"u}r $z$ \hbox{nahe 1,}
$$ \eq{y= f(z) &= (2r_0)^{1/2}(1-z)^{1/2}\(1+O(1-z)\)\,,\cr
f'(z) &=  - (\h r_0)^{1/2}(1-z)^{-1/2}\(1+O(1-z)\)\,,\cr f''(z) &=
- \h(\h r_0)^{1/2}(1-z)^{-3/2}\(1+O(1-z)\)\,.\cr} \eqno(4.0)  $$
Daraus folgt, da\ss\ $f'^2(z)+f(z)f''(z)$ auf ganz $[-1,1]$ stetig
ist, und daher
$$ M:=\max_{-1\le z\le 1}\abs{f'^2(z)+f(z)f''(z)}<\infty \eqno(4.1)  $$ gilt.
Weiters ersieht man f{\"u}r $z\to1-$, da\ss $f'(z)\to-\infty$, \
$f''(z)\to-\infty$ und \hbox{$f(z)f'(z)\to -r_0$.} Au\ss erdem folgt
aus dem Gesagten, f{\"u}r alle $z\in]-1,1[$,
$$ f(z)=f(-z)>0\,,\ \sgn f'(z)= -\sgn z\,,\ f''(z)<0\,,\qquad f(1)=0\,,\
f(0)\le1\,. \eqno(4.2)$$ Wir interessieren uns wieder f{\"u}r die
Gitterpunktanzahl im linear vergr{\"o}\ss erten K{\"o}rper $t\,\K$,
$t$ ein gro\ss er Parameter, also f{\"u}r
$$ A_{\K}(t) = \#\(t\,\K\cap\Z^3\) =
\# \left \{ (n_1,n_2,m) \in {\Bbb{Z}}^3: n_1^2 + n_2^2 \le t^2f^2
\left ( \frac{m}{t} \right ) , \abs{m} \le t \right \}\,, $$ und
insbesondere f{\"u}r den Gitterrest $
P_{\K}(t):=A_{\K}(t)-t^3\vol(\K)$.

\bsk\msk

{\bf Satz 3.}\quad{\it Es bezeichne $z_0\ge0$ die kleinste Zahl,
so da\ss\ $f''$ auf $[z_0,1[$ monoton ist\fn{6}{\rm Offensichtlich
ist $z_0<1$, da analog (4.0) $|f'''(z)|\asymp(1-z)^{-5/2}$ f\"ur
$z$ nahe 1 gilt.}, und es sei $\disp F_2:=\max_{[0,z_0[}|f''|$
falls $z_0>0$, und $F_2:=0$ f\"ur $z_0=0$. Dann gilt unter den
formulierten Voraussetzungen an den K\"orper $\K$ f{\"u}r jedes
$t>0$
$$ \abs{P_{\K}(t)} \le C_1\, t^{3/2} + C_2\, t\log t+C_3\,t+C_4 t^{3/4}
+ C_5\,t^{1/2} + C_6 $$ mit
$$ \eq{C_1  & := 53 \rf{1/2}{3/4} +12 \rf{1/2}{1/4} +73\,\rmax^{3/4} + 13\,
\rmax^{1/4} \,, \cr C_2  & := 7315\,\rmax^{1/2} + 241\,  \,, \cr
C_3  & := 3658\,\rmax^{1/2}\log_+(\rmax) +
642\,\log\({\rmax\over\rmin}\) + 1.6(M+r_0)+281\,, \cr C_4  & :=
19\rf{3/4}{9/8}+34\,\rmax^{5/8}  \,, \cr C_5  & :=
1302\,\rmax^{1/4}  \,, \cr C_6  & := 5  \,. \cr } $$ Dabei
bezeichnen $r_{\max}, r_{\min}$ die Extremwerte des
Kr{\"u}mmungsradius $r_\theta$ auf $\c$.} \bsk\bsk

{\bf Vorbereitung der Absch{\"a}tzung. } Als Pr{\"a}zisierung des
$O$-Symbols bezeichne $\Theta$ irgendeine reelle Gr{\"o}\ss e, die
von den Parametern beliebig abh{\"a}ngen darf, aber jedenfalls
betraglich $\le1$ ist. Offensichtlich gilt
$$ A_\K(t) = \sum_{\abs{m}\le t}\(\sum_{n_1^2 + n_2^2
\le t^2f^2 ( \frac{m}{t})}\ 1\)\,. \eqno(4.3) $$ Dabei ist die
innere Summe der Zahl der Gitterpunkte in einer Kreisscheibe.
F{\"u}r eine solche findet man im wesentlichen bereits bei Fricker
[5], S.~42/43, die elementare Absch{\"a}tzung (mit
$\psi(w):=w-[w]-\h$, wobei $[w]$ die gr{\"o}\ss te ganze Zahl $\le w$
bezeichnet)
$$ \sum_{n_1^2 + n_2^2 \le X^2}\ 1 \ = \pi X^2 -
8\sum_{0<n\le X/\sqrt{2}}\psi\(\sqrt{X^2-n^2}\) + 5\Theta\,,
\eqno(4.4)  $$ wenn man das dort angegebene $O(1)$ explizit macht.
Durch Einsetzen von (4.4) in (4.3) erh{\"a}lt man
$$ A_\K(t) = \pi t^2 \sum_{\abs{m}\le t} f^2\({m\over t}\)
- 8 \sum_{\abs{m}\le t}\(\sum_{0 < n \le
\frac{t}{\sqrt{2}}f(\frac{m}{t})} \psi \( \sqrt{t^2f^2 \(
\frac{m}{t} \) - n^2}\,\)\) + 5\Theta(2t+1)\,. $$ Auf die erste
Summe wird nun die Euler-MacLaurin'sche Summenformel angewendet.
(Vgl.~z.B.~Kr{\"a}tzel [\LP], S.~20.) Dies ergibt wegen
$f(\pm1)=0$
$$ \pi t^2 \sum_{\abs{m}\le t} f^2\({m\over t}\) =
\pi t^2 \Int_{-t}^t f^2 \left ( \frac{u}{t} \right ) \d u + 2\pi t
\Int_{-t}^t f \left( \frac{u}{t} \right) f' \left( \frac{u}{t}
\right) \psi (u) \d u\,. $$ Der erste Term rechts ist
offensichtlich $t^3\,\vol(\K)$. Mit $\psi_1(u) = \int_0^u \psi
(\tau )\d\tau = {1\over8}\Theta$  folgt wegen
$\disp\lim_{z\pm1}\,f(z)f'(z) = \mp r_0$ (vgl.~(4.1) - (4.2))

\vbox{$$ \Int_{-t}^t f \left( \frac{u}{t} \right) f' \left(
\frac{u}{t} \right) \psi (u) \d u ={\Theta\over4}\,r_0 - {1\over
t} \Int_{-t}^t \left(f'^2\ut +f\ut f''\ut\right) \psi_1(u) \d u =
$$
$$ = {\Theta\over4}\,r_0 - \Int_{-1}^1 \(f'^2(z)+f(z)f''(z)\)\psi_1(tz)\d z =
{\Theta\over4}\,(r_0+M)\,. $$ }

Wir erhalten also insgesamt
$$ A_\K(t) = \vol(\K)t^3-16\,P^*(t) + \(\(\h\pi(r_0+M)+10\)t+5\)\Theta\,,
\eqno(4.5)  $$
$$ P^*(t) := \sumi_{0\le m\le t}\(\sum_{0 < n \le \frac{t}{\sqrt{2}}
f(\frac{m}{t})}\psi \( \sqrt{t^2f^2 \( \frac{m}{t} \) -
n^2}\,\)\)\,, $$ wobei $\sum'$ bedeutet, da\ss\ der Summand mit $m=0$
den Faktor $\h$ erh{\"a}lt.

\bsk\msk

{\bf Hilfsmittel zum Beweis.}\msk

{\bf Hilfssatz 1. } (Exponentialsumme und $\psi$-Summe.) { \it
Durchl{\"a}uft $(n_1,n_2)$ eine beliebige endliche Teilmenge $J$
von $\Z^2$, dann gilt f\"ur jedes $Z>1$
$$ \eq{& \abs{ \sum_{(n_1,n_2) \in J} \psi (f(n_1,n_2)) }
\le  \cr & \le  \frac{1}{\pi z} \sum_{(n_1,n_2) \in J}1 +
\frac{1}{\pi} \sum_{\nu = 1}^{\infty} \min \left( \frac{1}{\nu},
\frac{Z^2}{\nu^3} \right) \abs{ \sum_{(n_1,n_2) \in J} e^{2\pi i
\nu f(n_1,n_2)}}\,.\cr} \leqno{{\rm(i)}}$$ Ferner bestehen die
Absch\"atzungen $$ \sum_{\nu = 1}^{\infty} \min \left(
\frac{1}{\nu}, \frac{Z^2}{\nu^3} \right)\,t^\alpha \le\ \cases{2Z
& f\"ur $\al=1$, \cr {8\over3}\,Z^{1/2} & f\"ur $\al=\h$, \cr \log
Z \, +2 & f\"ur $\al=0$.\cr  } \leqno{{\rm(ii)}}$$} \bsk

{\bf Beweis. }  Ohne numerische Konstanten findet man (i) in
Kr{\"a}tzel [\LP], Kap.~1.3, S.~26/27. Ein Beweis von (i) mit
numerischen Konstanten verl{\"a}uft v{\"o}llig analog zum Beweis
von Hilfssatz 1.3 in Kr{\"a}tzel [\AZ], S.~18. \par Um (ii) zu
beweisen, wendet man auf die Zerlegung $$ \sum_{\nu = 1}^{\infty}
\min \left( \frac{1}{\nu}, \frac{Z^2}{\nu^3} \right)\,t^\alpha = 1
+ \sum_{1<\nu\le Z}{1\over\nu} + Z^2 \sum_{\nu>Z} {1\over\nu^3} $$
zweimal die Eulersche Summenformel an (vgl.~[\LP], S.~20, Th.~1.3)
und sch\"atzt die auftretenden Restintegrale \"uber die Funktion
$\psi$ mit Hilfe des zweiten Mittelwertsatzes der Integralrechnung
ab. \qed

\bsk\bsk

{\bf Hilfssatz 2. } (Exponentialsumme und Exponentialintegral.) {
\it Es seien $a < b$ reell, $F$ stetig auf $[a,b]$, \ $F'$ stetig
in $]a, b[$ und monoton, $F'^*_{]a,b[}$ bezeichne das Bild von
$]a,b[$ unter $F'$. \ssk {\rm(i) } Falls
$F'^*_{]a,b[}\subseteq\,]0,1[$, folgt
$$ \sum_{a < n \le b}e^{2\pi iF(n)} = \int\limits_a^b e^{2\pi iF(u)}\d u +
\int\limits_a^b e^{2\pi i(F(u) - u)}\d u + \Theta\(1 +
\frac{4}{\pi}\)\,. $$ {\rm(ii) } Falls
$F'^*_{]a,b[}\subseteq\,]0,1-\varphi[$ f{\"u}r ein $\varphi>0$,
dann gilt
$$ \sum_{a < n \le b}e^{2\pi iF(n)} = \int\limits_a^b e^{2\pi iF(u)}\d u +
\Theta\,\(1 + \frac{2}{\pi} \left( \frac{1}{\varphi} + 1
\right)\)\,. $$ {\rm(iii) } Falls
$F'^*_{]a,b[}\subseteq\,]\varphi,1[$ f{\"u}r ein $\varphi>0$, dann
folgt
$$ \sum_{a < n \le b}e^{2\pi iF(n)} = \int\limits_a^b
e^{2\pi i(F(u) - u)}\d u + \Theta\,\(1 + \frac{2}{\pi}
\left(\frac{1}{\varphi} + 1 \right)\)\,. $$ } \bsk

{\bf Beweis.} (i) und (ii) werden in Wahrheit bereits in
Kr{\"a}tzel [\AZ], S.~167 hergeleitet, wenn man bei den Konstanten
nichts verschenkt. (iii) folgt aus (ii), indem man $F(u)$ durch $u
- F(u)$ ersetzt und konjugiert. \qed                    \bsk\msk

{\bf Hilfsatz 3. } { \it Es sei $F$ reell, zweimal stetig
differenzierbar auf $[a,b]$ und dort $|F''(u)| \ge \lambda > 0$.
Dann ist
$$ \left| \sum_{a < n \le b}e^{2\pi iF(n)} \right | \le
|F'(b) - F'(a)| \frac{5}{\sqrt{\lambda}} +
\frac{11}{\sqrt{\lambda}}\,. $$ } \bsk {\bf Beweis. } Das ist
(1.9) im Korollar zu Satz 1.3, S.~16 in Kr{\"a}tzel [\AZ].\qed

\bsk\bsk {\bf Hilfssatz 4.}\fn{7}{Dies ist eine effektive Version
der bereits in (4.0) enthaltenen Schranke
$f''(z)\ll(1-z)^{-3/2}$.} \quad{ \it Es sei $\c$ die Schnittkurve
des Randes $\partial\K$ mit der $(y,z)$-Ebene. Ihr Teil im ersten
Quadranten sei durch $y=f(z)\iff z=g(y)$ dargestellt. Dann gilt
f\"ur $0\le z<1$ $$ {1\over\rmax}\le |f''(z)| \le
8^{3/4}\rq(1-z)^{-3/2}\,.
$$ } \bsk

{\bf Beweis. } Aus $$
{1\over\rmax}\le{(1+f'^2(z))^{3/2}\over\rmax}\le |f''(z)| \le
{(1+f'^2(z))^{3/2}\over\rmin} \eqno(*) $$ folgt $$ |f''(z)|\le
\cases{{\sqrt{8}\over\rmin} & wenn $|f'(z)|\le1\,,$\cr
{\sqrt{8}\over\rmin}|f'(z)|^3 & wenn $|f'(z)|\ge1\,.$\cr }
\eqno(**)$$ F\"ur $|f'(z)|\le1$ ist damit mehr als die Behauptung
bewiesen. Wir setzen daher im Folgenden $|f'(z)|\ge1$ voraus, was
$|f'(u)|\ge1$ f\"ur $z\le u<1$ impliziert. Es sei $(y,z)\in\R_+^2$
mit $y=f(z)\iff z=g(y)$. Wegen $|f'(z)|\to\infty$ f\"ur $z\to1-$
ist
$$ {1\over f'(z)} = - \Int_1^z
{f''(u)\over f'^2(u)}\,\d u = - \Int_0^y {f''(g(w))\over
f'^3(g(w))}\,\d w\,, $$ daraus folgt mit $(*)$ \
$\disp{1\over|f'(z)|}\ge{y\over\rmax}$, also nach $(**)$
$$ |f''(z)|\le\sqrt{8}\,{\rmax^3\over\rmin}\,{1\over f^3(z)}\,. \eqno(***) $$
Wegen $g(0)=1$ ist weiters, mittels partieller Integration,
$$ 1-z=1-g(y)=-\Int_0^y g'(w)\d w = - \Int_0^y {\d w \over
f'(g(w))}= - \Int_0^y (y-w)\,{f''(g(w))\over f'^3(g(w))}\,\d w\,,
$$ folglich, wegen $(**)$, \ $1-z\le \sqrt{2}\,y^2\rmin^{-1}$,
also $y=f(z)\ge2^{-1/4}(\rmin\,(1-z))^{1/2}$. Somit ist nach
$(***)$ $$ |f''(z)|\le 8^{3/4}\,{\rmax^3\over
\rmin^{5/2}}\,(1-z)^{-3/2}\,.  $$ Da aus $g(0)=1$, $f(0)\le1$,
$f'(0)=g'(0)=0$ geometrisch $\rmin^{-1}\ge1$ folgt, ist damit
Hilfssatz 4 bewiesen. \qed

 \bsk\bsk

{\bf Absch\"atzung von Exponentialsummen. } Wegen (4.5) haben wir
Summen der Gestalt $$ T = T(a,b; t)  =  \sum_{at < m \le bt}\
\sum_{0 < n \le \frac{t}{\sqrt{2}} f(\frac{m}{t})} \psi \left (
\sqrt{t^2f^2 \left ( \frac{m}{t} \right )- n^2} \, \right)
\eqno(4.6)  $$ zu betrachten mit reellen Werten $0 < a < b \le 1$,
von denen die Absch\"atzung abh\"angen darf. Angesichts von
Hilfssatz 1 f\"uhrt dies auf Exponentialsummen
$$ \sum_{at < m \le bt}\
\sum_{0 < n \le \frac{t}{\sqrt{2}} f(\frac{m}{t})} e^{-2\pi i\nu
\sqrt{t^2f^2 \left ( \frac{m}{t} \right )- n^2}}\,, $$ $\nu$ eine
positive ganze Zahl. Wegen $$ 0 < n \le \frac{t}{\sqrt{2}}
f(\frac{m}{t})\quad\iff\quad 0 < {n\over\sqrt{x^2f^2(m/x) - n^2}}
\le 1 $$ erhalten wir in nat\"urlicher Weise Teilsummen
$$ S(k,\nu; a,b; t) := \sum_{at < m \le bt}\
\sum_{k<{\nu n/\sqrt{x^2f^2(m/x) - n^2}}\le k+1} e^{-2\pi i\nu
\sqrt{t^2f^2 \left ( \frac{m}{t} \right )- n^2}}\,,\eqno(4.7)
$$ $k\in\{0,1,\dots,\nu-1\}$. Wir spalten diese in weitere drei
Teilsummen $S_1, S_2, S_3$ auf, je nachdem, wie nahe bei der
jeweils inneren Summation ${\nu n/\sqrt{x^2f^2(m/x) - n^2}}$ der
n\"achst\-liegenden ganzen Zahl kommt. \bsk\bsk

{\bf Absch\"atzung der Exponentialsumme $S_1$. } F\"ur $k = 0,1,
\ldots , \nu -1$ und $0 < \vartheta \le \frac{1}{2}$ sei
$$ \eq{S_1 =S_1(k,\nu; a,b; t) &:=   \sum_{at < m \le bt}\
\sum_{k + \vartheta < \nu n/ \sqrt{x^2f^2(m/x) - n^2} \le k + 1 -
\vartheta} e^{- 2\pi i \nu \sqrt{t^2f^2(\frac{m}{t}) - n^2}} = \cr
& =  \sum_{at < m \le bt} \ \sum_{ \vartheta < \nu n/
\sqrt{t^2f^2(m/t) - n^2} - k \le 1 - \vartheta} e^{- 2\pi i( \nu
\sqrt{t^2f^2(\frac{m}{t}) - n^2} + kn)}\,.\cr}$$ Definieren wir
(bei durchwegs positiven reellen Werten aller Variabler) eine
Funktion $\varphi$ durch Aufl\"osen von
$${\d\over\d u_1}\(-\sqrt{t^2f^2(u_2/t) -u_1^2}\,\)= y\ \in]0,1]
$$ nach $u_1$, also $$ u_1={y\over\sqrt{1+y^2}}\,t\,f\({u_2\over t
}\)=:\varphi(y,u_2,t)\,,$$ dann folgt auch $$ {\d\over\d
u_1}\(-\nu \sqrt{t^2f^2(u_2/t) -u_1^2}\ -k u_1\,\)= y \ \iff\
u_1=\varphi\({k+y\over\nu},u_2,t\)\,.
$$ Somit ergibt sich $$
S_1 = \sum_{at < m \le bt}\ \sum_{\varphi ( \frac{k +
\vartheta}{\nu} , m,t) < n \le \varphi ( \frac{k + 1 -
\vartheta}{\nu}, m,t)} e^{- 2\pi i( \nu \sqrt{t^2f^2( \frac{m}{t})
- n^2} + kn)}\,.  $$ Nach Teil (i) von Hilfssatz 2 ist
$$ S_1 = S_{11} + S_{12} + \Theta\(1+{4\over\pi}\)((b-a)t+1) \eqno(4.8)$$
mit
$$ \eq{S_{11} & =  \sum_{at < m \le bt}\ \int\limits_{\varphi (
\frac{k + \vartheta}{\nu}, m,t)}^{ \varphi ( \frac{k + 1 -
\vartheta}{\nu}, m,t)} e^{- 2\pi i( \nu \sqrt{t^2f^2 (
\frac{m}{t}) - u^2} + ku)} \d u,  \cr S_{12} & =  \sum_{at < m \le
bt} \int\limits_{\varphi ( \frac{k + \vartheta}{\nu}, m,t)}^{
\varphi ( \frac{k + 1 - \vartheta}{\nu}, m,t)} e^{- 2\pi i( \nu
\sqrt{t^2f^2 ( \frac{m}{t}) - u^2} + (k + 1)u)} \d u\,.  \cr }
$$
Wir f\"uhren die Absch\"atzung von $S_{11}$ im Detail aus, jene
von $S_{12}$ verl\"auft genauso und liefert das gleiche Ergebnis.
Die Substitution
$$ u = \varphi \left ( \frac{k + y}{\nu}, m,t \right ) =
\frac{k + y}{\sqrt{\nu^2 + (k + y)^2}} \, t\,f \left ( \frac{m}{t}
\right )$$ (mit $y$ als neuer Integrationsvariabler) ergibt
$$ \nu \sqrt{t^2f^2 \left ( \frac{m}{t} \right ) - u^2}\ +ku =
\eta(y) \, t\,f \left ( \frac{m}{t} \right )\,,\quad   \eta(y) :=
\frac{\nu^2+k(k+y) }{\sqrt{\nu^2 + (k + y)^2}} $$ mit $$
\frac{\nu}{\sqrt{2}} \le \frac{\nu^2}{\sqrt{\nu^2 + (k + 1)^2}}
\le \eta(y)\le \sqrt{\nu^2+(k+y)^2} \le \sqrt{2} \, \nu\,,
\eqno(4.9)$$ und weiter nach einfacher Rechnung $\d u = -t
f\(\frac{m}{t}\)\,{\eta'(y)\over y}\,\d y$, also
$$ \int\limits_{\varphi (
\frac{k + \vartheta}{\nu}, m,t)}^{ \varphi ( \frac{k + 1 -
\vartheta}{\nu}, m,t)} e^{- 2\pi i( \nu \sqrt{t^2f^2 (
\frac{m}{t}) - u^2} + ku)} \d u = -t\,f\(\frac{m}{t}\)
\Int_\vartheta^{1-\vartheta}e^{-2\pi i\eta(y)t
f(m/t)}\,{\eta'(y)\over y}\,\d y\,. $$ Mittels partieller
Integration ist dies gleich
$$ -\frac{1}{2\pi i} \, \frac{1}{1 -
\vartheta}\,  e^{-2 \pi i\eta(1 - \vartheta )t\,f( \frac{m}{t})} +
\frac{1}{2\pi i} \, \frac{1}{\vartheta}\,  e^{-2 \pi i\eta(
\vartheta )t\,f(\frac{m}{t})} - \frac{1}{2\pi i}
\Int_\vartheta^{1-\vartheta}e^{-2\pi i\eta(y)t f(m/t)}\,{\d y\over
y^2}\,.$$ Damit erhalten wir insgesamt
$$ S_{11} = - S_{111} + S_{112} - S_{113} \eqno(4.10)$$ mit
$$ \eq{S_{111} & =  \frac{1}{2\pi i} \, \frac{1}{1 -
\vartheta} \sum_{at < m \le bt} e^{-2 \pi i\eta(1 - \vartheta
)t\,f( \frac{m}{t})}, \cr S_{112} & =  \frac{1}{2\pi i} \,
\frac{1}{ \vartheta} \sum_{at < m \le bt} e^{-2 \pi i\eta(
\vartheta )t\,f( \frac{m}{t})}, \cr  S_{113} & = \frac{1}{2\pi i}
\Int_\vartheta^{1-\vartheta} \(\sum_{at < m \le bt} e^{-2\pi
i\eta(y)t f(m/t)}\)\,{\d y\over y^2}\,.\cr}$$ Auf jede dieser
Exponentialsummen wird jetzt Hilfssatz 3 angewendet, mit
$$ F(u) = -\eta(y)\,t\,f \left ( \frac{u}{t} \right )\,, \quad \quad
y\ \in\ [\vartheta , 1 - \vartheta]\,. $$ Dann ist f\"ur $at\le
u\le bt$
$$ \eq{F'(u) &= -\eta(y)f'\left( \frac{u}{t} \right )>0\,,\cr
F''(u) &= -\eta(y)\,{1\over t}\,f''\left( \frac{u}{t} \right ) \ge
\eta(y)\,{1\over t}\,\f =:\la>0\,.\cr} $$

Mit $0 < \vartheta \le \frac{1}{2}$ folgt nach Hilfssatz 3, unter
Beachtung von (4.9),
$$ \eq{|S_{111}| & \le  \frac{5}{\pi} \, \eta(1 - \vartheta )
\abs{f'(b) - f'(a)} \sqrt{ \frac{t}{\eta(1 - \vartheta )\f }} +
\cr
  & + \frac{11}{\pi} \sqrt{ \frac{t}{\eta(1 - \vartheta )\f }}\le  \cr
& \le  \frac{5}{\pi} \, 2^{\frac{1}{4}}|f'(b) - f'(a)| \sqrt{
\frac{\nu t}{\f }} + \frac{11}{\pi} \, 2^{\frac{1}{4}} \sqrt{
\frac{t}{\nu \f }}\,,   \cr |S_{112}| & \le  \frac{5}{\pi} \, 2^{
- \frac{3}{4}} \, \frac{1} { \vartheta} \, |f'(b) - f'(a)| \sqrt{
\frac{\nu t}{\f }} + \frac{11}{\pi} \, 2^{ - \frac{3}{4}} \,
\frac{1}{\vartheta} \, \sqrt{ \frac{t}{\nu \f }}\,. \cr |S_{113}|
& \le \( \frac{5}{\pi} \, 2^{- \frac{3}{4}} |f'(b) - f'(a)| \sqrt{
\frac{\nu t}{\f }} + \frac{11}{\pi} \, 2^{- \frac{3}{4}} \sqrt{
\frac{t}{\nu \f }}\, \) \frac{1}{\vartheta}\,.\cr}  $$ Somit
ergibt sich f\"ur $S_{11}$ nach (4.10) insgesamt
$$ |S_{11}| \le \frac{5}{\pi} \, 2^{\frac{1}{4}} |f'(b) -
f'(a)| \sqrt{ \frac{\nu t}{\f }} \left \{ \frac{1}{\vartheta} + 1
\right \} + \frac{11}{\pi} \, 2^{\frac{1}{4}} \sqrt{ \frac{t}{\nu
\f }} \left \{ \frac{1}{\vartheta} + 1 \right \}.     $$ Da f\"ur
$S_{12}$ die gleiche Absch\"atzung gilt, folgt mit (4.9)
$$ \eq{|S_1| & \le  \frac{10}{\pi} \, 2^{\frac{1}{4}} |f'(b) -
f'(a)| \sqrt{ \frac{\nu t}{\f }} \left \{ \frac{1}{\vartheta} + 1
\right \} +  \cr &  + \frac{22}{\pi} \, 2^{\frac{1}{4}} \sqrt{
\frac{t}{\nu \f }} \left \{ \frac{1}{\vartheta} + 1 \right \} +
\left ( 1 + \frac{4}{\pi} \right )((b - a)t + 1)\,.\cr }
\eqno(4.11)  $$ \bsk\bsk

{\bf Absch\"atzung der Exponentialsumme $S_2$. } F\"ur $0 <
\vartheta \le \frac{1}{2}$ und $k = 0,1,\ldots , \nu -1$
sei\fn{8}{Ein Blick zur\"uck auf (4.7) ist hier hilfreich.}
$$ S_2=S_2(k,\nu; a,b; t) =
\sum_{at < m \le bt} \sum_{k < \nu n/ \sqrt{t^2f^2 (m/t) - n^2}
\le k + \vartheta} e^{-2\pi i \nu \sqrt{t^2f^2( \frac{m}{t}) -
n^2}}\,. $$ Nach Teil (ii) von Hilfssatz 2 ist
$$ S_2 = S_{21} + \Theta\,\left ( 1 + \frac{6}{\pi} \right )((b - a)t + 1)
\eqno(4.12)  $$ mit
$$ S_{21} :=  \sum_{at < m \le bt}\ \int\limits_{\varphi (
\frac{k}{\nu}, m,t)}^{\varphi (\frac{k + \vartheta} {\nu}, m,t)}
e^{-2\pi i( \nu \sqrt{t^2f^2( \frac{m}{t}) - u^2} + ku)}\d u \,.
$$ Die Substitution
$$ u = \varphi \left ( \frac{k + y}{\nu}, m,t \right ) =
\frac{k + y}{\sqrt{\nu^2 + (k + y)^2}} \, tf \left ( \frac{m}{t}
\right ) $$ im Integral ergibt wie vorher
$$ S_{21}  = -  \sum_{at < m \le bt}\
\int\limits_{0}^{\vartheta}  tf \left ( \frac{m}{t} \right )
e^{-2\pi i\eta(y)tf( \frac{m}{t})}\, {\eta'(y)\over y}\,\d y = $$
$$  = - \int\limits_{0}^{\vartheta} {\eta'(y)\over y}\,
\left \{ tf(b) \sum_{at < m \le bt} e^{-2\pi i\eta(y)tf(
\frac{m}{t})} -
  \int\limits_{at}^{bt} f' \left ( \frac{u}{t} \right )
\sum_{at < m \le u} e^{-2\pi i\eta(y)tf( \frac{m}{t})}\d u \right
\} \d y\,.$$ mittels partieller Summation. \msk Auf die
Exponentialsummen hier wird nun wieder Hilfssatz 3 angewendet.
Beachtet man $f(b) \le f(a) \le f(0) \le 1$, so ergibt sich
$$ \eq{|S_{21}| & \le \abs{\int\limits_{0}^{\vartheta} {\eta'(y)\over y}\,
\d y}\, \left \{ 10 \cdot 2^{\frac{1}{4}} |f'(b) - f'(a)|
\sqrt{\frac{\nu t^3}{\f }} + 22 \cdot
2^{\frac{1}{4}}\sqrt{\frac{t^3}{\nu \f }}\, \right \}  \cr & \le
\vartheta \left \{ 10 \cdot 2^{\frac{1}{4}} |f'(b) - f'(a)|
\sqrt{\frac{t^3}{\nu \f }} + 22 \cdot 2^{\frac{1}{4}}
\sqrt{\frac{t^3}{\nu^3 \f }} \, \right \}\,.\cr }   $$ Zusammen
mit (4.12) erhalten wir daher
$$ \eq{|S_2| & \le  \vartheta \left \{ 10 \cdot 2^{\frac{1}{4}} |f'(b) - f'(a)|
\sqrt{\frac{t^3}{\nu \f }} + 22 \cdot 2^{\frac{1}{4}}
\sqrt{\frac{t^3}{\nu^3 \f }} \, \right \} +  \cr &  + \left ( 1 +
\frac{6}{\pi} \right )((b - a)t + 1)\,. \cr } \eqno(4.13) $$
\bsk\bsk

{\bf Absch\"atzung der Exponentialsumme $S_3$. } Wir betrachten
zuletzt
$$ S_3 = \sum_{at <m \le bt} \sum_{k + 1 - \vartheta < \nu n/
\sqrt{t^2f^2(m/t) - n^2} \le k + 1} e^{-2\pi i \nu \sqrt{t^2f^2(
\frac{m}{t}) - n^2}} $$ mit $0 < \vartheta \le \frac{1}{2}$,
$k=0,1,\dots,\nu-1$. Nach Teil (iii) von Hilfssatz 2 erh\"alt man
analog der eben durchgef\"uhrten Absch\"atzung
$$ \eq{|S_3| & \le \vartheta \left \{ 10 \cdot 2^{\frac{1}{4}} |f'(b) - f'(a)|
\sqrt{\frac{t^3}{\nu \f }} + 22 \cdot 2^{\frac{1}{4}}
\sqrt{\frac{t^3}{\nu^3 \f }} \, \right \} +  \cr &  + \left ( 1 +
\frac{6}{\pi} \right )((b - a)t + 1)\,. \cr }  \eqno(4.14) $$
\bsk\bsk

{\bf Zusammenfassung der Exponentialsummen-Absch\"atzung. } Durch
Kombination der Ergebnisse (4.11), (4.13) und (4.14) erhalten wir
f\"ur die in (4.7) erkl\"arte Summe
$$ \eq{ \abs{S(k,\nu; a,b; t)} & \le
\frac{10}{\pi} \cdot 2^{\frac{1}{4}}|f'(b) - f'(a)|
\sqrt{\frac{\nu t}{\f }} \left \{ \frac{1}{\vartheta} + 1 \right
\} +  \cr &  + \frac{22}{\pi} \cdot 2^{\frac{1}{4}}
\sqrt{\frac{t}{\nu \f }} \left \{ \frac{1}{\vartheta} + 1 \right
\} +  \cr } $$
$$  \eq{+\ \vartheta & \( 20 \cdot 2^{\frac{1}{4}}|f'(b) - f'(a)|
\sqrt{\frac{t^3}{\nu \f }} + 44 \cdot 2^{\frac{1}{4}}
\sqrt{\frac{t^3}{\nu^3 \f }} \, \) +\cr    &+ \left ( 3 +
\frac{16}{\pi} \right )((b - a)t + 1)\,.\cr}   $$

W\"ahlt man nun f\"ur $\nu \le \frac{\pi t}{2}$
$$ \vartheta = \sqrt{\frac{\nu}{2\pi t}}\,, $$
so ergibt sich
$$ \eq{|S(k,\nu; a,b; t)| & \le \frac{20}{\sqrt{\pi \f }} \, 8^{\frac{1}{4}}
|f'(b) - f'(a)|t + \cr & + \frac{10}{\pi} \, 2^{\frac{1}{4}}|f'(b)
- f'(a)| \sqrt{\frac{\nu t}{\f }} +  \cr &  + \frac{44}{\sqrt{\pi
\f }} \, 8^{\frac{1}{4}} \cdot \frac{t}{\nu} + \frac{22}{\pi} \,
2^{\frac{1}{4}} \sqrt{\frac{t}{\nu \f }} + \cr & + \left ( 3 +
\frac{16}{\pi} \right )((b - a)t + 1)\,, \cr } $$ was sich leicht
zu
$$ \eq{|S(k,\nu; a,b; t)| & \le \frac{20}{\sqrt{\pi \f }} \, 8^{\frac{1}{4}}
|f'(b) - f'(a)|t + \cr & + \frac{10}{\pi} \, 2^{\frac{1}{4}}|f'(b)
- f'(a)| \sqrt{\frac{\nu t}{\f }} +  \cr &  + \frac{55}{\sqrt{\pi
\f }} \, 8^{\frac{1}{4}} \cdot \frac{t}{\nu} + \left ( 3 +
\frac{16}{\pi} \right )((b - a)t + 1)\cr } \eqno(4.15)  $$
vereinfacht. Die Richtigkeit dieser letzten Absch\"atzung f\"ur
$\nu > \frac{\pi t}{2}$ folgt sofort allein aus den Schranken
f\"ur $S_2$ und $S_3$ mit $\vartheta = \frac{1}{2}$. \bsk\bsk

{\bf Absch\"atzung der $\psi$-Summe. } Nach diesen Vorbereitungen
sind wir in der Lage, die vor (4.6) erkl\"arte Bruchteilsumme
$T(a,b;t)$ zu behandeln. Angesichts von (4.6), (4.7) und (4.14),
verwendet f\"ur $k=0,1,\dots,\nu-1$, folgt mittels Hilfssatz 1,
mit noch verf\"ugbaren Parametern $\nu\in\Z^+$ und $Z>1$,
$$ \eq{|T(a,b;t)| & \le  \frac{1}{\pi Z} \sum_{at < m \le bt}\
\sum_{0 <n\le \frac{t}{\sqrt{2}}}\ 1\ +  \cr  &  + \frac{1}{\pi}
\sum_{\nu = 1}^{\infty} \min \left ( \frac{1}{\nu},
\frac{Z^2}{\nu^3} \right ) \left | \sum_{at < m \le bt}\ \sum_{0 <
n \le \frac{t}{\sqrt{2}} f( \frac{m}{t})} e^{- 2\pi i \nu
\sqrt{t^2f^2( \frac{m}{t}) - n^2}} \right |  \cr } $$ $$ \eq{& \le
\frac{1}{\pi \sqrt{2}} ((b - a)t + 1) \, \frac{t}{Z} +  \cr &  +
\frac{1}{\pi} \sum_{\nu = 1}^{\infty} \min \left ( \frac{1}{\nu},
\frac{Z^2}{\nu^3} \right ) \left \{ \frac{20 \nu} {\sqrt{\pi \f }}
\, 8^{\frac{1}{4}} |f'(b) - f'(a)| t +  \right.  \cr &   +
\frac{10}{\pi} \, 2^{\frac{1}{4}} |f'(b) - f'(a)| \sqrt{\frac{\nu
t}{\f }} + \frac{55}{\sqrt{\pi \f }} \, 8^{\frac{1}{4}}t +  \cr &
 \left. + \left ( 3 + \frac{16}{\pi} \right ) ((b - a)t + 1)
\nu  \right \}  \cr  & \le \frac{b - a}{\pi \sqrt{2}} \,
\frac{t^2}{Z} + \frac{40Z} {\sqrt{\pi^3 \f }} \, 8^{\frac{1}{4}}
|f'(b) - f'(a)|t +  \cr &   + \frac{80}{3\pi^2} \, 2^{\frac{1}{4}}
|f'(b) - f'(a)| \sqrt{\frac{Z t}{\f }} + \frac{55}{\sqrt{\pi^3 \f
}} \, 8^{\frac{1}{4}} t( \log Z\, + 2) + \cr & + \frac{1}{\pi}
\left ( 6 + \frac{32}{\pi} \right ) ((b - a)t + 1)Z + \frac{t}{\pi
\sqrt{2}}\,.\cr} $$ Wir w\"ahlen nun
$$ Z := \frac{1}{4\sqrt{5}} \, 2^{- \frac{1}{8}}(\pi \f
)^{\frac{1}{4}} \sqrt{\frac{(b - a)t}{|f'(b) - f'(a)|}}\,.
\eqno(4.16) $$ (F\"ur den Fall, da\ss\ $Z\le1$ vgl.~man (4.18)
unten.) Dann erhalten wir aus der zuletzt angeschriebenen
Absch\"atzung
$$ \eq{
|T(a,b;t)| & \le 4\sqrt{5} \left ( \frac{2}{\pi^2} \right
)^{\frac{5}{8}} \sqrt{(b - a) |f'(b) - f'(a)|} \, \ff ^{-
\frac{1}{4}} t^{\frac{3}{2}} +  \cr &  + \frac{40}{3} \,
2^{\frac{3}{16}} \, 5^{- \frac{1}{4}} \pi^{- \frac{15}{8}} (b -
a)^{\frac{1}{4}} |f'(b) - f'(a)|^{\frac{3}{4}} \ff ^{-
\frac{3}{8}} t^{\frac{3}{4}} + \cr &  + \frac{55}{\sqrt{\pi^3 \f
}}\, 2^{- \frac{1}{4}}t \log \left ( \frac{(b - a)\sqrt{\f } \,
t}{|f'(b) -f'(a)|} \right ) +  \cr  & + \frac{1}{2\sqrt{5}} \,
2^{- \frac{1}{8}} \pi^{- \frac{3}{4}} \left ( 3 + \frac{16}{\pi}
\right ) ((b - a)t + 1) \ff ^{\frac{1}{4}} \sqrt{\frac{(b -
a)t}{|f'(b) - f'(a)|}}  \cr & + \frac{t}{\pi \sqrt{2}}\,.\cr }  $$
Um dies weiter zu vereinfachen, ben\"utzen wir die folgenden
Ungleichungen:
$$ \eq{&
\sqrt{(b - a)|f'(b) - f'(a)|} \, \ff ^{- \frac{1}{4}} = \sqrt{(b -
a) \int\limits_a^b |f''(u)|\d u} \, \ff ^{- \frac{1}{4}} \le  \cr
& \le (b - a) \sqrt{\F } \, \ff ^{- \frac{1}{4}}  \le \left |
\frac{\F }{\f } \right |^{1/2} \int\limits_a^b
|f''(u)|^{\frac{1}{4}}\d u\,,  \cr}$$ $$ \eq{ (b -
a)^{\frac{1}{4}}|f'(b) - f'(a)|^{\frac{3}{4}} \ff ^{- \frac{3}{8}}
& \le \left | \frac{\F }{\f } \right |^{3/4} \int\limits_a^b
|f''(u)|^{\frac{3}{8}}\d u, \cr (b - a) \sqrt{\frac{b - a}{|f'(b)
- f'(a)|}} \, \ff ^{\frac{1}{4}} & \le \left | \frac{\F }{\f }
\right |^{1/4} \int\limits_a^b |f''(u)|^{- \frac{1}{4}}\d u\,.\cr
} $$ Weiters bemerken wir: Entweder ist $$ \frac{(b - a) \sqrt{\f
}}{|f'(b) - f'(a)|} \le 1\,, $$ dann ist der entsprechende
Logarithmus kleiner oder gleich $0$, oder es gilt $$  1 \le
\frac{(b - a) \sqrt{\f }}{|f'(b) - f'(a)|} \le \frac{1}{\sqrt{\f
}} \le \sqrt{r_{\max}}\,. $$ Damit erhalten wir insgesamt
$$ \eq{& 16|T(a,b;t)|  \le \cr & \le \left ( 53 \left |
\frac{\F }{\f } \right |^{\frac{1}{2}} \int\limits_a^b
|f''(u)|^{\frac{1}{4}}\d u + 12 \left|\frac{\F }{\f }
\right|^{\frac{1}{4}} \int\limits_a^b |f''(u)|^{- \frac{1}{4}}\d u
\right ) t^{\frac{3}{2}} + \cr & + 19 \left | \frac{\F }{\f }
\right |^{\frac{3}{4}}\,t^{3\over4}\, \int\limits_a^b
|f''(u)|^{\frac{3}{8}}\d u  + \frac{133}{\sqrt{\f }} \, t \left (
\log t + \frac{1}{2} \log_+(r_{\max}) \right ) +  \cr & + 12\ff
^{- \frac{1}{4}} \sqrt{t} + 6t\,,\cr}\eqno(4.17)$$ mit $\log_+ :=
\max(0,\log)$. \ssk Sollte $Z\le 1$ sein, dann ergibt sich mittels
trivialer Absch\"atzung
$$\eq{16|T(a,b;t)| & \le \frac{8}{\sqrt{2}}((b - a)t^2 + t) \cr
& \le \frac{8}{\sqrt{2}} \, 4 \sqrt{5} \, 2^{\frac{1}{8}} \pi^{-
\frac{1}{4}} \sqrt{(b - a) |f'(b) - f'(a)|} \ff ^{-
\frac{1}{4}}t^{\frac{3}{2}} + 6t  \cr & \le 42 \left | \frac{\F
}{\f } \right |^{\frac{1}{2}} \int\limits_a^b
|f''(u)|^{\frac{1}{4}}\d u + 6t\,.\cr} \eqno(4.18) $$ Damit hat
das obige Ergebnis (4.17) Allgemeing\"ultigkeit.

\bsk\bsk

{\bf Vollendung der Restabsch\"atzung. } Wir setzen bei (4.5) fort
und sch\"atzen die Summanden f\"ur $m=0$ und $m\in]t-1,t]$ trivial
ab. Ihr Beitrag zu $P^*(t)$ ist
$\Theta\,{\text{3\over8}}\sqrt{2}\,t$, daher gilt
$$ A_\K(t) = \vol(\K)t^3-16\,P^{**}(t)
+ \(\(\h\pi(r_0+M)+6\sqrt{2}+10\)t+5\)\Theta\,,   $$
$$ P^{**}(t) := \sum_{0< m\le t-1}\(\sum_{0 < n \le \frac{t}{\sqrt{2}}
f(\frac{m}{t})}\psi \( \sqrt{t^2f^2 \( \frac{m}{t} \) -
n^2}\,\)\)\,. \eqno(4.19) $$

Es bezeichne $z_0\ge0$ die kleinste Zahl, so da\ss\ $f''$ auf
$[z_0,1[$ monoton ist. Wir setzen weiters
 $c:=\({54\over53}\)^2$ und definieren
eine Folge $(z_j)_{1\le j\le J}$ durch die Gleichung $|f''(z_j)| =
c^j |f''(z_0)|$. Dabei sei $J$ maximal, so da\ss\ $z_J t\le t- 1$,
also $1-z_J\ge t^{-1}$. Wegen $\rmax^{-1}\,c^J\le
|f''(z_J)|\le8^{3/4}\rq(1-z_J)^{-3/2}$ (nach Hilfssatz 4) folgt
mit kurzer Rechnung
$$ J\le {3\,\log t\over2\log c}+ C\,,\qquad
C:={1\over\log c}\,\log\(8^{3/4}\,{\rmax^4\over\rmin^4}\)\,.
\eqno(4.20) $$ Nach unserer Konstruktion ist klarerweise
$$ P^{**}(t) = T(0,z_0; t) + \sum_{j=0}^J T(z_j, z_{j+1}; t)\,. \eqno(4.21) $$
Wir sch\"atzen zuerst $T(0,z_0; t)$ ab. Dazu sei $\disp
F_2:=\max_{[0,z_0[}|f''(z)|$ definiert. Wir verwenden (4.17) mit
$a=0, b=z_0$, sowie die Schranken $\rmax^{-1}\le|f''(z)|\le F_2$
und erhalten
$$ \eq{16&|T(0,z_0; t)| \le \(53\,\rf{1/2}{3/4}+12\,\rf{1/2}{1/4}\)\,t^{3/2}
+ 133\,\rmax^{1/2}\,t \log t + \cr & +
(66.5\,\rmax^{1/2}\,\log_+(\rmax) + 6)\,t +
19\,\rf{3/4}{9/8}\,t^{3/4} + 12\,\rmax^{1/4}\,t^{1/2}\,.\cr }
\eqno(4.22)  $$ Zur Absch\"atzung von $T(z_j, z_{j+1}; t)$ wenden
wir (4.17) mit $a=z_j$, $b=z_{j+1}$ an, dann ist $ {\F\over\f}=c$,
$\f=c^j$, und mit Hilfssatz 4 ergibt sich
$$ \eq{16&|T(z_j, z_{j+1} ; t)| \le \cr \le &\
\(53\,c^{1/2} \iz{1/4} \d w +12\,c^{1/4}\,\rmax^{1/4}
(z_{j+1}-z_j) \)\,t^{3/2} + \cr & + 133\,c^{-j/2}\,\rmax^{1/2}\,t
\log t + \(66.5\,c^{-j/2}\,\rmax^{1/2}\,\log_+(\rmax) + 6\)\,t +
\cr & + 19\,c^{3/4}\,t^{3/4} \,\iz{3/8} \d w +
12\,\rmax^{1/4}\,c^{-j/4}\,t^{1/2}\,.\cr }  \eqno(4.23)  $$ Beim
Aufsummieren \"uber $j=0,1,\dots,J$ beachten wir, da\ss
$$ 133\,\sj c^{-j/2} < 7182\, \quad\und\quad 12\,\sj c^{-j/4}<1290 $$ gilt,
und erhalten
$$ \eq{16&\sj |T(z_j, z_{j+1} ; t)| \le \cr \le &\
\(53\,c^{1/2} \Int_0^1 |f''(w)|^{1/4} \d w +
12\,c^{1/4}\,\rmax^{1/4} \)\,t^{3/2} + \cr & +
7182\,\rmax^{1/2}\,t \log t + \(3591\,\rmax^{1/2}\,\log_+(\rmax) +
6(J+1)\)\,t + \cr & + 19\,c^{3/4}\,t^{3/4} \,\Int_0^1
|f''(w)|^{3/8} \d w + 1290\,\rmax^{1/4}\,t^{1/2}\,.\cr }
\eqno(4.24)  $$ Die auftretenden Integrale werden durch den
maximalen Kr\"ummungsradius abgesch\"atzt: Mit $\be={1\over4},\
{3\over8}$ ist
$$ \Int_0^1 |f''(w)|^\be \d w \le \Int_0^1 |f''(w)|\,{\rmax^{1-\be}
\over(1+f'^2(w))^{3(1-\be)/2}}\,\d w = $$ $$ =
\rmax^{1-\be}\,\Int_0^\infty {\d u\over(1+u^2)^{{3(1-\be)/2}}} =
\rmax^{1-\be}\, {\sqrt{\pi}\,\Gamma(1-{3\over2}\be)\over2\,
\Gamma({3\over2}(1-\be))} \approx \cases{1.35\,\rmax^{3/4} & f\"ur
$\be={1\over4}$,\cr 1.725\,\rmax^{5/8} & f\"ur
$\be={3\over8}$.\cr}  $$ Wir m{\"u}ssen nur mehr die Resultate
(4.19), (4.21), (4.22) und (4.24) zusammenf{\"u}gen, f{\"u}r $J$
die Schranke (4.20) sowie den Wert $c:=\({54\over53}\)^2$
einsetzen und vereinfachen. Dies ergibt nach etwas m{\"u}hsamer
Rechnung (z.B.~unterst{\"u}tzt von {\it Derive } [15]) die in Satz
3 dargestellte Absch{\"a}tzung des Gitterrestes. \qed \bsk {\bf
Bemerkung. } F\"ur den Spezialfall der Kugel ergibt Satz 3 (wegen
$F_2=0$, $\rmax=1$) als gr\"o\ss ten Fehlerterm $86\,t^{3/2}$.
Dies ist nat\"urlich etwas schw\"acher als $14\,t^{3/2}$ von Satz
1A, aber es zeigt, da\ss\ die wesentlich verschiedene und viel
allgemeinere Beweismethode des Satzes 3 bez\"uglich der Sch\"arfe
in etwa derselben Gr\"o\ss enordnung bleibt.

\vbox{\vskip 2.5true cm}

\klein \parindent=0pt

\cen{\bf Literatur}  \bsk \def\smc{}

[1] {\smc V.~Bentkus \and F.~G\"otze, } On the lattice point
problem for ellipsoids. Acta Arithm. {\bf80}, 101--125 (1997).  \ssk

[2] {\smc S.~Bochner, } Die Poisson'sche Summenformel in mehreren
Ver\"anderlichen. Math.~Ann. {\bf106}, 56--63 (1932). \ssk

[3] {\smc F.~Chamizo, } Lattice points in bodies of revolution.
Acta Arith. {\bf85}, 265-277 (1998). \ssk

[4] {\smc F.~Chamizo \and H.~Iwaniec, } On the sphere problem.
Rev.~Mat.~Iberoamericana {\bf 11}, 417-429 (1995). \ssk

[5] {\smc F.~Fricker, } Einf\"uhrung in die Gitterpunktlehre.
Basel-Boston-Stuttgart 1982. \ssk

[6] {\smc F.~G\"otze, } Lattice point problems and values of quadratic forms, 
Preprint, Univ.~Bielefeld 2004.  \ssk 

[7] {\smc I.S.~Gradshteyn \and I.M.~Ryzhik, } Table of integrals,
series, and products. A.~Jeffrey editor. 5th ed., San Diego 1994.
\ssk

[8] {\smc D.R.~Heath-Brown, } Lattice points in the sphere. In:
Number theory in progress, Proc. Number Theory Conf. Zakopane
1997, eds. K.~Gy\"ory et al., vol.~{\bf2}, 883-892 (1999).  \ssk

[9] {\smc E.~Hlawka, } \"Uber Integrale auf konvexen K\"orpern I.
Monatsh.~f.~Math. {\bf 54}, 1-36 (1950); II, ibid.~{\bf 54},
81--99 (1950). \ssk

[10] M.~Huxley, Exponential sums and lattice points III.
Proc.~London Math.~Soc., III.~Ser., {\bf87} (2003), 591-609. \ssk 

[11] {\smc E.~Kr\"atzel, } Lattice points. Berlin 1988. \ssk

[12] {\smc E.~Kr\"atzel, } Analytische Funktionen in der
Zahlentheorie. Stuttgart-Leipzig-Wiesbaden 2000.  \ssk

[13] {\smc E.~Kr\"atzel, } Lattice points in convex planar
domains. Monatsh.~Math., im Druck.  \ssk

[14] {\smc W.~M\"uller, } Lattice points in large convex bodies.
Monatsh.~Math. {\bf128}, 315-330 (1999).  \ssk

[15] {\smc Soft Warehouse, } Inc.,  {\it Derive,} Version 3.11,
Honolulu (Hawaii) 1995.  \ssk

[16] {\smc J.G.~Van der Corput, } Zahlentheoretische
Absch\"atzungen mit Anwendungen auf Gitterpunktprobleme. Math.Z.
{\bf 17}, 250--259 (1923).  \ssk

[17] {\smc I.M.~Vinogradov, } On the number of integer points in a
sphere (Russian). Izv.~Akad.~Nauk SSSR Ser.~Mat. {\bf27}, 957-968
(1963). \ssk

[18] {\smc A.~Walfisz, } Weylsche Exponentialsummen in der neueren
Zahlentheorie. Berlin 1963.

\vbox{\vskip 1.5true cm}

\parindent=1.5true cm

\vbox{Ekkehard Kr\"atzel \ssk

Institut f\"ur Mathematik

Universit\"at Wien

Nordbergstra\ss e 15 

1090 Wien, \"Osterreich \ssk

http://www.univie.ac.at/\~baxa/kraetzel.html

\bsk\msk

Werner Georg Nowak \ssk

Institut f\"ur Mathematik

Department f\"ur Integrative Biologie

Universit\"at f\"ur Bodenkultur Wien

Peter Jordan-Stra\ss e 82

1190 Wien, \"Osterreich \ssk

E-mail: {\tt \ nowak@mail.boku.ac.at} \ssk

http://www.boku.ac.at/math/nth.html}

\bye